\documentclass[11pt]{article}
\usepackage[dvips, hmargin=1.5in, vmargin={1.6in, 1.6in}]{geometry}
\usepackage{amssymb, amsmath, amsthm, stmaryrd, enumerate}
\usepackage[ps,dvips,arrow, matrix,tips,curve]{xy}
\SelectTips{cm}{10}

\DeclareMathOperator{\ob}{ob}
\newcommand{\db}[1]{\mathbb{#1}}
\newcommand{\cat}[1]{\mathbf{#1}}
\newcommand{\op}{\mathrm{op}}

\newcommand{\id}{\mathrm{id}}
\newcommand{\thg}{{\mathord{\text{--}}}}

\newcommand{\abs}[1]{{\left|{#1}\right|}}
\newcommand{\dbr}[1]{\left[\,{#1}\,\right]}

\newcommand{\spn}[1]{{\left<{#1}\right>}}
\newcommand{\elt}[1]{\left\ulcorner{#1}\right\urcorner}

\newcommand{\defn}[1]{\textbf{#1}}
\newcommand{\cd}[2][]{\vcenter{\hbox{\xymatrix#1{#2}}}}
\newcommand{\cdl}[2][]{\xymatrix@1#1{#2}}

\renewcommand{\b}[1]{\boldsymbol{\mathbf{#1}}}
\newcommand{\I}{{\b I}}

\newcommand{\C}{{\b C}}
\newcommand{\D}{{\b D}}
\newcommand{\M}{\mathbb M}
\renewcommand{\P}{\mathbb P}

\newcommand{\X}{{\b X}}
\newcommand{\Y}{{\b Y}}

\newcommand{\xtor}[1]{\cdl[@1]{{} \ar[r]|-{\object@{|}}^{#1} & {}}}
\newcommand{\tor}{\ensuremath{\relbar\joinrel\mapstochar\joinrel\rightarrow}}
\newcommand{\To}{\ensuremath{\Rightarrow}}
\newcommand{\Tor}{\ensuremath{\Relbar\joinrel\Mapstochar\joinrel\Rightarrow}}

\newcommand{\twocong}[2][0.5]{\ar@{}[#2] \save ?(#1)*{\cong}\restore}
\newcommand{\rtwocell}[3][0.5]{\ar@{}[#2] \ar@{=>}?(#1)+/l 0.2cm/;?(#1)+/r 0.2cm/^{#3}}
\newcommand{\ltwocell}[3][0.5]{\ar@{}[#2] \ar@{=>}?(#1)+/r 0.2cm/;?(#1)+/l 0.2cm/_{#3}}
\newcommand{\dtwocell}[3][0.5]{\ar@{}[#2] \ar@{=>}?(#1)+/u  0.2cm/;?(#1)+/d 0.2cm/^{#3}}
\newcommand{\dthreecell}[3][0.5]{\ar@{}[#2] \ar@3{->}?(#1)+/u  0.2cm/;?(#1)+/d 0.2cm/^{#3}}
\newcommand{\utwocell}[3][0.5]{\ar@{}[#2] \ar@{=>}?(#1)+/d 0.2cm/;?(#1)+/u 0.2cm/_{#3}}
\newcommand{\dtwocelltarg}[3][0.5]{\ar@{}#2 \ar@{=>}?(#1)+/u  0.2cm/;?(#1)+/d 0.2cm/^{#3}}
\newcommand{\utwocelltarg}[3][0.5]{\ar@{}#2 \ar@{=>}?(#1)+/d  0.2cm/;?(#1)+/u 0.2cm/_{#3}}
\newcommand{\set}[2]{\left\{\,#1 \mid #2\,\right\}}

\theoremstyle{plain}
\newtheorem{Prop}{Proposition}
\newtheorem{Cor}[Prop]{Corollary}
\newtheorem{Thm}[Prop]{Theorem}

\theoremstyle{definition}
\newtheorem{Defn}[Prop]{Definition}
\newtheorem{Rk}[Prop]{Remark}

\newcommand{\coll}{{\db Coll}}

\newcommand{\dCat}{{\db Cat}}
\newcommand{\Cat}{{\cat{Cat}}}

 \makeatletter \@namedef{itemize*}{\itemize\parsep\z@ \parskip\z@}
 \@namedef{enditemize*}{\enditemize} \@namedef{enumerate*}{\enumerate\parsep\z@
 \parskip\z@} \@namedef{endenumerate*}{\endenumerate} \makeatother

\begin{document}
\title{Polycategories via pseudo-distributive laws}
\author{Richard Garner\thanks{Supported by a PhD grant from the EPSRC}\\
Department of Pure Mathematics, University of Cambridge,\\Wilberforce Road,
Cambridge CB3 0WB, UK} \maketitle

\begin{abstract}
In this paper, we give a novel abstract description of Szabo's
\emph{polycategories}. We use the theory of \emph{double
  clubs}~--~a generalisation of Kelly's theory of clubs to `pseudo'
(or `weak') double categories~--~to construct a pseudo-distributive law of the
free symmetric strict monoidal category pseudocomonad on $\cat{Mod}$ over
itself \emph{qua} pseudomonad, and show that monads in the `two-sided Kleisli
bicategory' of this pseudo-distributive law are precisely symmetric
polycategories.
\end{abstract}


\section{Introduction}

Szabo's theory of \emph{polycategories} \cite{szabo:poly} has been the target
of renewed interest over recent years. Polycategories are the
`not-necessarily-representable' cousins of the \emph{weakly
  distributive categories} of \cite{seely:wdc}; their relationship
mirrors that of \emph{multicategories} to \emph{monoidal categories}.

Though it is possible, as Szabo did, to give a `hands on' description of a
polycategory, such a description leaves a lot to be desired. For a start, the
sheer quantity of data that one must check for even simple proofs quickly
becomes overwhelming. Further problems arise when one wishes to address aspects
of a putative `theory of polycategories': what are the correct notions of
polyfunctor or polytransformation? What is a polycategorical limit?  In
attempting to answer such questions without a formal framework, one is forced
into the unsatisfactory position of relying on intuition alone.

Thus far, the paper \cite{koslowski:poly} has provided the only attempt to
rectify this situation. Koslowski provides an abstract description of
polycategories that generalises the elegant work of \cite{burroni:tcat} and
later \cite{hermida:rep-multi} and \cite{leinster:higheroperads} on
`$T$-multicategories'. However, whilst this latter theory uses only some rather
simple and obvious constructions on categories with finite limits, the
structures that Koslowski uses to build his description of polycategories are
rather more complicated and non-canonical. Furthermore, the generalisation from
the non-symmetric to the symmetric case is not as smooth as one would like.

We therefore offer an alternative approach to the abstract description of
polycategories.  It is the same and not the same as Koslowski's: again, we
shall build on an abstract description of multicategories, and again,
composition proceeds using something like a `distributive law'. Where we
deviate from Koslowski is in the description of multicategories that we build
upon.

In Section 1, we recount this alternative description: it is the approach of
\cite{baezdolan:ncat} and \cite{tholen:multicats}, based on \emph{profunctors}
rather than \emph{spans}. We go on to describe how we may generalise this
description to one for polycategories; to do this we invoke a
\emph{pseudo-distributive law} (in the sense of \cite{marm:psd1},
\cite{tanaka:thesis}) of a pseudocomonad (the `target arity') over a
pseudomonad (the `source arity'). Polycategories now arise as monads in the
`two-sided Kleisli bicategory' of this pseudo-distributive law.

There are several advantages to this approach: it allows us to describe
\emph{symmetric} polycategories with no greater difficulty than
\emph{non-symmetric} polycategories; it will generalise easily from ordinary
categories to enriched categories; and, though we do not attempt this here, it
allows us to `read off' further aspects of the theory of polycategories: the
aforementioned polyfunctor, polytransformation, and so on.

In order to make this description go through, we must construct a suitable
pseudo-distributive law. Now, a pseudo-distributive law is a prodigiously
complicated object: it is five pieces of (complex) data subject to ten
coherence laws. A bare hands construction would be both tedious and
unenlightening: the genuinely interesting combinatorics involved would be
obscured by a morass of trivial details.

Thus, in Section 2, we discuss how we may use the theory of \emph{double
clubs}, as developed in the companion paper \cite{rhgg1}, to reduce this
Herculean task to something more manageable. Informally, the theory of double
clubs tells us that it suffices to construct our pseudo-distributive law at the
terminal category $1$, and that we can propagate this construction elsewhere by
`labelling objects and arrows' appropriately.

Finally, in Section 3, we perform this construction at $1$; and though one
might think this would be an exercise in nose-following, it actually turns out
to be a fairly interesting piece of categorical combinatorics. Equipped with
this, we are finally able to prove the existence of our pseudo-distributive law
and hence to give our preferred definition of polycategory.

An Appendix gives the definitions of \emph{pseudomonad}, \emph{pseudocomonad}
and \emph{pseudo-distributive law}.

\section{Multicategories and polycategories}\label{chap:poly}

We begin by re-examining the theory of multicategories: the material here
summarises \cite{baezdolan:ncat}, \cite{hyland:prooftheory} and
\cite{leinster:operads}, amongst others. Note that throughout, we shall only be
interested in the theory of \emph{symmetric} multicategories, and, later, of
\emph{symmetric} polycategories: that is, we allow ourselves to reorder freely
the inputs and outputs of our maps. Consequently, whenever we say
`multicategory' or `polycategory', it may be taken that we mean the symmetric
kind. The non-symmetric case for polycategories is considered in more detail by
\cite{koslowski:poly}.

\subsection{Multicategories}
We write $X^*$ for the free monoid on a set $X$, and $\Gamma, \Delta, \Sigma,
\Lambda$ for typical elements thereof. We will use commas to denote the
concatenation operation on $X^*$, as in ``$\Gamma, \Delta$''; and we will tend
to conflate elements of $X$ with their image in $X^*$. Given $\Gamma = x_1,
\dots, x_n \in X^*$, we define $\abs{\Gamma} = n$, and given $\sigma \in S_n$,
write $\sigma \Gamma$ for the element $x_{\sigma(1)}, \dots, x_{\sigma(n)} \in
X^*$.

\begin{Defn}
A \defn{symmetric multicategory} $\M$ consists of:
\begin{itemize*}
\item A set $\ob \M$ of \defn{objects};
\item For every $\Gamma \in (\ob \M)^*$ and $y \in \ob \M$, a set $\M(\Gamma;
    y)$ of \defn{multimaps} from $\Gamma$ to $y$ (we write a typical element of
    such as $f \colon \Gamma \to y$); further, for every $\sigma \in
    S_{\abs{\Gamma}}$, an
\defn{exchange isomorphism}
$\M(\Gamma; y) \to \M(\sigma \Gamma; y)$.
\item For every $x \in \ob \M$, an \defn{identity map} $\id_x \in \M(x; x)$;
\item For every $\Gamma, \Delta_1, \Delta_2 \in (\ob \M)^*$ and $y, z \in \ob
    \M$, a \defn{composition map}
\[\M(\Gamma; y) \times \M(\Delta_1, y, \Delta_2; z) \to \M(\Delta_1, \Gamma, \Delta_2; z)\text,\]
\end{itemize*}
This data satisfies axioms expressing the fact that exchange isomorphisms
compose as expected, and that composition is associative, unital, and
compatible with exchange isomorphisms: see \cite{lambek:multi} for the full
details.
\end{Defn}
\noindent  Now, this data expresses composition as a \emph{binary} operation
performed between two multimaps; however, there is another view, where we
`multicompose' a family of multimaps $g_i \colon \Gamma_i \to y_i$ with a
multimap $f \colon y_1, \dots, y_n \to z$.

The transit from one view to the other is straightforward: we recover the
multicomposition from the binary composition by performing, in any order, the
binary compositions of the $g_i$'s with $f$: the axioms for binary composition
ensure that this gives a uniquely defined composite. Conversely, we can recover
binary composition from multicomposition by setting all but one of the $g_i$'s
to be the identity.
\def\itimes{,}

We can express the operation of multicomposition as follows: fix the object set
$X = \ob \M$, and consider it as a discrete category. We write $S$ for the
\emph{free symmetric strict monoidal category} 2-monad on $\cat{Cat}$, and
consider the functor category $[(SX)^\op \times X, \cat{Set}]$. To give an
object $F$ of this is to give sets of multimaps as above, together with
coherent exchange isomorphisms. Further, this category has a `substitution'
monoidal structure given by
\[(G \otimes F)(\Gamma; z) = \!\!\!\!\!\!\!\sum_{\substack{k \in \mathbb N\\y_1, \dots, y_k \in X}}\!\!\!\!\!\int^{\Delta_1, \dots, \Delta_k \in SX} \!\!\!\!\!\!\!\!\!\!G(y_1, \dots, y_k; z) \times \prod_{i=1}^k F(\Delta_i; y_i) \times SX(\Gamma, \bigotimes_{i=1}^k \Delta_i)\text,\]
and
\[\I(\Gamma; x) = \begin{cases}\{\ast\} & \text{if $\Gamma = x$}\\\emptyset &
\text{otherwise;}\end{cases}\]

and to give a multicategory is precisely to give a monoid with respect to this
monoidal structure. Indeed, suppose we have a monoid $F \in [(SX)^\op \times X,
\cat{Set}]$. Then the unit map $j \colon \I \to F$ picks out for each $x \in X$
an element of $F(x; x)$, which will correspond to the identity multimap $\id_x
\colon x \to x$. What about the multiplication map $m \colon F \otimes F \to
F$? Unpacking the above definition, we see that $(F \otimes F)(\Gamma; z)$ can
be described as follows. Let $\Delta_1, \dots, \Delta_k \in (\ob \M)^*$ be such
that
\begin{itemize*}
\item $\abs{\Gamma} = n = \sum \abs{\Delta_i}$;
\item there exists $\sigma \in S_n$ such that $\sigma \Gamma = \Delta_1, \dots,
    \Delta_k$,
\end{itemize*}
and let $f_i \colon \Delta_i \to y_i$ (for $i = 1, \dots, k$), and $g \colon
y_1, \dots, y_k \to z$ be multimaps in $F$. Then this gives us a typical
element of $(F \otimes F)(\Gamma; z)$, which we visualise as
\[\cd[@C=0em]{\Gamma \ar[d]^{\sigma} \\
 \Delta_1, \dots, \Delta_k \ar[d]^{f_1, \dots, f_k}\\
 y_1 \itimes \dots \itimes y_k \ar[d]^g\\
 z\text.
}
\]

The map $m \colon F \otimes F \to F$ sends this element to an element of
$F(\Gamma; z)$; in other words, it specifies the result of this
`multicomposition'. The associativity and unitality laws for a monoid ensure
that this composition process is associative and unital as required.

In fact, we may deduce the existence of the substitution monoidal structure on
$[(SX)^\op \times X, \cat{Set}]$ from more abstract considerations. The key
idea is to construct a bicategory $\mathcal B$ with $\mathcal B(X, X) =
[(SX)^\op \times X, \cat{Set}]$, in such a way that horizontal composition in
this endohom-category induces the desired substitution monoidal structure; and
for this, we make use of the following result:

\begin{Prop}\label{lifting}
The symmetric strict monoidal category 2-monad $(S, \eta, \mu)$ on $\cat{Cat}$
lifts to a pseudomonad $(\hat S, \hat \eta, \hat \mu, \lambda, \rho, \tau)$ on
$\cat{Mod}$, the bicategory of categories, profunctors and transformations.
\end{Prop}
(For the definition of and notation for a pseudomonad, see the Appendix).
\begin{proof} We recount only the salient details here. For a full proof the
reader may refer to \cite{tanaka:thesis}; but see also Section
\ref{presentation} below.

The lifted homomorphism $\hat S \colon \cat{Mod} \to \cat{Mod}$ agrees with $S
\colon \cat{Cat} \to \cat{Cat}$ on objects; whilst on 1-cells, it sends the
profunctor $F \colon \D^\op \times \C \to \cat{Set}$ to the profunctor $\hat S
F \colon (S \D)^\op \times S\C \to \cat{Set}$ given by:
\begin{equation*}
    \hat S F\big((d_1, \dots, d_n), (c_1, \dots, c_m)\big) =
    \begin{cases}
    \displaystyle\sum_{\sigma \in S_n} \prod_{i = 1}^n F(d_i, c_{\sigma(i)}) & \text{if $n = m$;} \\
    \qquad \quad 0 & \text{otherwise.}    \end{cases}
\end{equation*}
The components at $\C$ of the lifted transformations $\hat \eta \colon
\id_\cat{Mod} \Rightarrow \hat S$ and $\hat \mu \colon \hat S\hat S \Rightarrow
\hat S$ are obtained as the images of the corresponding components of $\eta$
and $\mu$ under the canonical embedding $(\thg)_\ast \colon \cat{Cat} \to
\cat{Mod}$. Explicitly, we have:
\begin{align*}
    \hat \eta_\mathbf C \big( \Gamma, c\big) & = S\C\big(\Gamma, (c)\big)\text;\\
    \hat \mu_\mathbf C \big( \Gamma, (\Delta_1, \dots, \Delta_n)\big) &=
    S\C(\Gamma, \bigotimes \Delta_i)\text.
\end{align*}\qedhere
\end{proof}

%

Now, just as each monad on a category gives rises to a Kleisli category, so
each pseudomonad on a bicategory gives rise to a `Kleisli bicategory'. This
construction was first given in \cite{hylandpowercheng:pseudo} for the special
case of a pseudomonad on a 2-category; and the following is the obvious
generalisation to the bicategorical case:
\begin{Defn}\label{kleislidef}
Let $\mathcal B$ be a bicategory and let $(S, \eta, \mu, \lambda, \rho, \tau)$
be a pseudomonad on $\mathcal B$. Then the Kleisli bicategory $Kl(S)$ of the
pseudomonad $S$ has:
\begin{itemize*}
\item \textbf{Objects} those of $\mathcal B$;
\item \textbf{Hom-categories} given by $Kl(S)(X, Y) = \mathcal B(X, SY)$;
\item \textbf{Identity map} at $X$ given by the component $\eta_X \colon X \to
    S X$;
\item \textbf{Composition} $Kl(S)(Y, Z) \times Kl(S)(X, Y) \to Kl(S)(X, Z)$
    given by
\[\cd{\mathcal B (Y, S Z) \times \mathcal B (X, S Y)
  \ar[d]^{\cong} \\
 1 \times \mathcal B (Y, S Z) \times \mathcal B (X, S Y)
  \ar[d]^{\elt{\mu_Z} \times S \times \id}\\
 \mathcal B (S S Z, S Z) \times \mathcal B (S Y, S S Z) \times \mathcal B (X, S Y)
  \ar[d]^{\otimes}\\
 \mathcal B (X, S Z)}
\]
where we use $\otimes$ to stand for some choice of order of composition for
this threefold composite. Explicitly, on maps, this composition is given by
\[(Y \xrightarrow{G} SZ)  \otimes (X \xrightarrow{F} SY)
 =
 X \xrightarrow{F} SY \xrightarrow{S G} SSZ \xrightarrow{\mu_Z} S Z\]
for some choice of bracketing for this composite.
\end{itemize*}
\end{Defn}
\noindent The remaining data to make this a bicategory~--~namely, the
associativity and unitality constraints~--~can be constructed in an obvious way
using the associativity and unitality constraints for $\mathcal B$ and the
coherence modifications for the pseudomonad $S$. The reader may easily verify
that these data satisfy the bicategory axioms.

\begin{Rk}
We may justify the name `Kleisli bicategory' as follows. At the one-dimensional
level, the Kleisli category of a monad $S$ on a category $\C$ is determined by
its universality amongst all categories $\D$ equipped with an embedding functor
$H \colon \C \to \D$ and a right action $\theta \colon HS \Rightarrow H$ of $S$
on $H$. Similarly, we may characterise the `Kleisli bicategory' of a
pseudomonad $S$ on a bicategory $\mathcal B$ as universal amongst all
bicategories $\mathcal D$ equipped with an embedding pseudofunctor $H \colon
\mathcal B \to \mathcal D$ and a right pseudo-action $\theta \colon HS
\Rightarrow S$ of $S$ on $H$: see \cite{hylandpowercheng:pseudo}, Theorem 4.3.
\end{Rk}

In particular, we may form the Kleisli bicategory of the pseudomonad $\hat S$
on $\cat{Mod}$; and by substituting the data given in the proof of Proposition
\ref{lifting} into Definition \ref{kleislidef}, we may easily verify that
horizontal composition in $Kl(\hat S)(X, X)$ gives precisely the monoidal
structure on $[(SX)^\op \times X, \cat{Set}]$ described above. Hence we arrive
at an alternative, but equivalent, definition of multicategory:
\begin{Defn}
A symmetric multicategory is a monad on a discrete object $X$ in the bicategory
$Kl(\hat S)$.
\end{Defn}
\noindent This description is well known, though not often stated in precisely
this form: it is the approach of \cite{baezdolan:ncat} and
\cite{tholen:multicats}.

\subsection{Polycategories}\label{polycatssec}
We recall now the notion of symmetric \emph{polycategory}:
\begin{Defn}\label{defpoly}
A \defn{symmetric polycategory} $\P$ consists of
\begin{itemize}
\item A set $\ob \P$ of \defn{objects};
\item For each pair $(\Gamma, \Delta)$ of elements of $(\ob \P)^\ast$, a
    set $\P(\Gamma; \Delta)$ of \defn{polymaps} from $\Gamma$ to $\Delta$;
\item For each $\Gamma$, $\Delta \in (\ob \P)^*$, each $\sigma \in
    S_{\abs{\Gamma}}$ and $\tau \in S_{\abs{\Delta}}$, \defn{exchange
    isomorphisms}
\[\P(\Gamma; \Delta) \to \P(\sigma \Gamma; \tau \Delta)\text,\]
\item For each $x \in \ob \P$, an \defn{identity map} $\id_x \in \P(x; x)$;
\item For $\Gamma, \Delta_1, \Delta_2, \Lambda_1, \Lambda_2, \Sigma \in
    (\ob \P)^\ast$, and $x \in \ob \P$, \defn{composition maps}
\[\P(\Gamma; \Delta_1, x, \Delta_2) \times \P(\Lambda_1, x, \Lambda_2; \Sigma) \to \P(\Lambda_1, \Gamma, \Lambda_2; \Delta_1, \Sigma, \Delta_2)\text,\]
\end{itemize}
subject to laws expressing the associativity and unitality of composition,
expressing that the exchange isomorphisms compose as expected, and that they
are compatible with composition: see \cite{szabo:poly} or \cite{seely:wdc} for
the full details.
\end{Defn}
\noindent We recover the notion of a multicategory if we assert that
$\P(\Gamma; \Delta)$ is empty unless $\Delta$ is a singleton.
\def\outimes{\circ}

Now, as before, we may shift from giving a `binary composition' of two polymaps
to giving a `polycomposition' operation on two families of composable polymaps.
First, we need to say what we mean by \emph{composable}.
\begin{Defn}
Let $\b{f} := \{f_m \colon \Lambda_m \to \Sigma_m\}_{1 \leqslant m \leqslant
j}$ and $\b{g} := \{g_n \colon \Gamma_n \to \Delta_n\}_{1 \leqslant n \leqslant
k}$ be families of polymaps,
such that
\[\sum \abs{\Sigma_m} = \sum \abs{\Gamma_n} = l\text.\]
We say that a permutation $\sigma \in S_l$ is a \defn{matching} of $\b f$ and
$\b g$ if $\sigma(\Sigma_1, \dots, \Sigma_j) = \Gamma_1, \dots, \Gamma_k$.
\end{Defn}
Informally, a matching of two families $\b f$ and $\b g$ indicates `which
output of $f_i$ has been plugged into which input of $g_j$'. Yet not every such
plugging need be obtainable from repeated binary composition; and so if our
notion of polycomposition is to have the same force as our notion of binary
composition, we must restrict the matchings along which we will allow
polycomposition to occur.

\begin{Defn}
Given a matching $\sigma$ of $\b{f}$ and $\b{g}$, we define a bipartite
multigraph $G_\sigma$ as follows. Its two vertex sets are labelled by $f_1,
\dots, f_m$ and $g_1, \dots, g_n$, and we add one edge between $f_i$ and $g_j$
for every element of $\Sigma_i$ which is paired with an element of $\Gamma_j$
under the matching $\sigma$.
%
%
We shall say that the matching $\sigma$ is
\defn{suitable} just when $G_\sigma$ is acyclic, connected and has no multiple
edges.
\end{Defn}

\begin{Prop}\label{full}
Let there be given families $\b f := \{f_m \colon \Lambda_m \to \Sigma_m\}_{1
\leqslant m \leqslant j}$ and $\b{g} := \{g_n \colon \Gamma_n \to \Delta_n\}_{1
\leqslant n \leqslant k}$ of polymaps; together with a suitable matching
$\sigma$ thereof. Then there is a uniquely defined polymap $\b g \circ_\sigma
\b f \colon \Lambda_1, \dots, \Lambda_j \to \Delta_1, \dots, \Delta_k$ obtained
by repeated binary compositions which, in some order, connect each $x \in
\Sigma_1, \dots, \Sigma_j$ with the corresponding $\sigma(x) \in \Gamma_1,
\dots, \Gamma_k$.
\end{Prop}


To prove this, we will prove something slightly stronger. First, a little more
notation: given a list $\Sigma = x_1, \dots, x_k \in X^\ast$, by a
\emph{sublist} of $\Sigma$ we shall mean a list $\Gamma = x_{i_1}, \dots,
x_{i_j}$ where $1 \leqslant i_1 < i_2 < \dots < i_j \leqslant k$. Thus sublists
of $\Sigma$ are in bijection with subsets of $\{1, \dots, \abs{\Sigma}\}$, and
in particular, form a Boolean algebra; and we write $\Gamma^c$ for the
complement of $\Gamma$ in this Boolean algebra. We also say that a list
$\Sigma$ is an \emph{interleaving} of two lists $\Gamma_1$ and $\Gamma_2$ if we
can view $\Gamma_1$ and $\Gamma_2$ as complementary sublists of $\Sigma$.

\begin{Defn}
Let $\b{f} := \{f_m \colon \Lambda_m \to \Sigma_m\}_{1 \leqslant m \leqslant
j}$ and $\b{g} := \{g_n \colon \Gamma_n \to \Delta_n\}_{1 \leqslant n \leqslant
k}$ be families of polymaps. A \defn{partial matching} of $\b f$ and $\b g$ is
given by a sublist $\Sigma$ of $\Sigma_1, \dots, \Sigma_m$ and a sublist
$\Gamma$ of $\Gamma_1, \dots, \Gamma_n$ with $\abs{\Sigma} = \abs{\Gamma} = l$,
together with a permutation $\sigma \in S_l$ satisfying $\sigma(\Sigma) =
\Gamma$.
\end{Defn}
As before, we can define the notion of the associated graph $G_\sigma$ for a
partial matching, and thus the notion of a \emph{suitable} partial matching.
Proposition \ref{full} now follows \emph{a fortiori} from the following:
\begin{Prop}
Let there be given families of polymaps $\b f$ and $\b{g}$ as before, together
with a suitable partial matching $(\Sigma, \Gamma, \sigma)$ thereof. Then there
is a uniquely defined polymap $\b g \circ_\sigma \b f$ obtained by repeated
binary compositions which, in some order, connect each $x \in \Sigma$ to the
corresponding $\sigma(x) \in \Gamma$. The domain of $\b g \circ_\sigma \b f$ is
an interleaving of the lists $\Lambda_1, \dots, \Lambda_j$ and $\Gamma^c$,
whilst the codomain is an interleaving of the lists $\Delta_1, \dots, \Delta_k$
and $\Sigma^c$.
\end{Prop}

\begin{proof}
Since the partial matching $\sigma$ is suitable, its associated graph
$G_\sigma$ is a tree, and so in particular will have a vertex of degree $1$.
Choose any such vertex: it corresponds to one of our polymaps $f_i$ or $g_i$,
without loss of generality to $f_i$, say. We begin by forming the binary
composition of $f_i$ with the polymap $g_j$ which is connected to $f_i$ in
$G_\sigma$. Suppose
\[f_i \colon \Lambda_i \to \Sigma_i, x, \Sigma'_i \quad \text{and}
\quad g_j \colon \Gamma_j, x, \Gamma'_j \to \Delta_j\text,\] where
the two $x$'s are matched under $\sigma$. Then the resultant composite map will
be
\[g_i \circ f_j \colon \Gamma_j, \Lambda_i, \Gamma'_j \to
\Sigma_i, \Delta_j, \Sigma_i'\text.\] Note that $f_i$ has no other outputs
taking part in the partial matching $\sigma$. Thus we can now form a partial
matching $\sigma'$ of $\b f \setminus \{f_i\}$ with $\b g \setminus \{g_j\}
\cup \{g_j \circ f_i\}$, which simply matches elements in the same way as
$\sigma$ except for the no-longer present matching of $x$. Now it's easy to see
that the associated graph of $\sigma'$ will be the same as that of $\sigma$,
but with the vertex corresponding to $f_i$ and the single adjacent edge
removed. We continue by induction on the size of the tree $G_\sigma$.

Note that we may at each stage have several possible choices of vertices of
degree $1$ which we may take as the next binary composition to perform.
However, the associativity laws for a polycategory ensure that the resultant
composite will be independent of the choice we make at each stage.  \qedhere
\end{proof}

Thus, in any polycategory, we may define the `polycomposition' of a family
$\b{f}$ with a family $\b{g}$ along a suitable matching $\sigma$: conversely,
if we are given polycomposites along suitable matchings, we may recapture a
binary composition by polycomposing with a suitable collection of identity
maps. Consequently, if we are to give an abstract formulation of polycategory,
it seems reasonable to do so in terms of a notion of `polycompositional'
polycategory.

In order to fully justify this last claim, we must exhibit a bijection between
polycompositional polycategories and the polycategories of Definition
\ref{defpoly}. However, we do not yet have a full description of the axioms
which a polycompositional polycategory should satisfy; and to write them down
at this point would be very messy. Thus we postpone justification until we have
given our abstract description of polycompositional polycategories, from which
we will be able to extract a description of the axioms such a structure must
satisfy; and hence to prove that these entities coincide with the
polycategories of Definition~\ref{defpoly}.

To arrive at our abstract formulation, we imitate the methods of the previous
section. Firstly, given a set $X$ of objects, we may view it as a discrete
category and consider the functor category $[(SX)^\op \times SX, \cat{Set}]$;
and to give an element of this is to give sets of polymaps together with
coherent exchange isomorphisms. We would now like to set up a monoidal
structure on this category such that a monoid in it is precisely a
polycompositional polycategory. The unit is straightforward:
\[I(\Gamma; \Delta) = \begin{cases}\{\ast\} & \text{if $\Gamma = x = \Delta$}\\\emptyset &
\text{otherwise;}\end{cases}\] and we can describe what a typical
element of $(F \otimes F)(\Gamma; \Delta)$ should look like. Let
\[
\Psi_1, \dots, \Psi_k\ \text, \quad \Lambda_1, \dots, \Lambda_k\ \text, \quad
\Sigma_1, \dots, \Sigma_l \quad \text{and} \quad \Phi_1, \dots, \Phi_l \] be
elements of $(\ob \M)^\ast$, such that
\begin{itemize*}
\item $\abs{\Gamma} = n = \sum \abs{\Psi_i}$;
\item $\sum \abs{\Lambda_i} = m = \sum \abs{\Sigma_j}$;
\item $\sum \abs{\Phi_j} = p = \abs{\Gamma}$;
\item there exists $\sigma \in S_n$ such that $\sigma \Gamma = \Psi_1, \dots,
    \Psi_k$;
\item there exists $\tau \in S_m$ such that $\tau$ is a suitable matching of
    $\{\Lambda_i\}$ with $\{\Sigma_i\}$;
\item there exists $\upsilon \in S_p$ such that $\upsilon(\Phi_1, \dots,
    \Phi_k) = \Delta$;
\end{itemize*}
and let $f_i \colon \Psi_i \to \Lambda_i$ (for $i = 1, \dots, k$), and $g_j
\colon \Sigma_j \to \Phi_j$ (for $j = 1, \dots, l$) be polymaps in $F$. Then
this gives us a typical element of $(F \otimes F)(\Gamma; \Delta)$, which we
visualise as
\[\cd[@C=0em]{\Gamma \ar[d]^{\sigma} \\
 \Psi_1, \dots, \Psi_k \ar[d]^{f_1, \dots, f_k}\\
 \Lambda_1, \dots, \Lambda_k \ar[d]^{\tau}\\
 \Sigma_1, \dots, \Sigma_l \ar[d]^{g_1, \dots, g_l}\\
 \Phi_1, \dots, \Phi_l \ar[d]^{\upsilon}\\
 \Delta\text.
}
\]

Then as for the multicategory case, the multiplication map $m \colon F \otimes
F \to F$ should specify a composite map for this `formal polycomposite', and
the associativity and unitality conditions for a monoid should ensure that this
polycomposition is associative and unital.

So our problem is reduced to finding a suitable way of expressing this monoidal
structure; and in fact we will skip straight over this stage and instead
describe polycompositional polycategories as monads in a suitable bicategory.
For this, we shall need the following fact:
\begin{Prop}\label{liftco}
The 2-monad $(S, \eta, \mu)$ on $\cat{Cat}$ lifts to a pseudocomonad $(\hat S,
\hat \epsilon, \hat \Delta)$ as well as a pseudomonad $(\hat S, \hat \eta, \hat
\mu)$ on $\cat{Mod}$.
\end{Prop}
\begin{proof}
The transformations $\hat \epsilon$ and $\hat \Delta$ have respective
components at $\cat{C}$ given by
\[{\hat \epsilon}_{\cat C} = (\eta_\cat{C})^\ast \quad \text{and} \quad
{\hat \Delta}_{\cat C} = (\mu_\cat{C})^\ast\text.\] We obtain the remaining
data for the pseudocomonad via the calculus of mates \cite{kellystreet:review},
making use of the adjunctions ${\hat \eta}_{\cat{C}} \dashv {\hat
\epsilon}_{\cat{C}}$ and ${\hat \mu}_{\cat{C}} \dashv {\hat \Delta}_{\cat{C}}$.
\qedhere
\end{proof}

[Since the pseudomonad $(\hat S, \hat \eta, \hat \mu)$ and the pseudocomonad
$(\hat S, \hat \epsilon, \hat \mu)$ share the same underlying homomorphism
$\hat S \colon \cat{Mod} \to \cat{Mod}$, there is some scope for confusion
here. To remedy this, we will use $\hat S_m$ and $\hat S_c$ as aliases for the
homomorphism $\hat S$; the former when we are thinking of it as part of a
pseudomonad structure, and the latter, when as part of a pseudocomonad.]

The key idea is to produce a \emph{pseudo-distributive law} $(\delta, \overline
\eta, \overline \epsilon, \overline \mu, \overline \Delta)$ of the
pseudocomonad $\hat S_c$ over the pseudomonad $\hat S_m$; that is, there should
be a pseudo-natural transformation $\delta \colon \hat S_c\hat S_m \Rightarrow
\hat S_m\hat S_c$ satisfying the rules of a distributive law `up to
isomorphism', as specified by the  invertible modifications $\overline \eta$,
$\overline \epsilon$, $\overline \mu$ and $\overline \Delta$: for full details,
see the Appendix.
Given such a pseudo-distributive law, polycategories will emerge as monads in
its `two-sided Kleisli bicategory'. Since this construction may not be
familiar, we describe it first one dimension down:
\begin{Defn}
Let $\C$ be a category, let $(S, \eta, \mu)$ be a monad and $(T, \epsilon,
\Delta)$ a comonad on $\C$, and let $\delta \colon TS \Rightarrow ST$ be a
distributive law of the comonad over the monad; so we have the four equalities:
\begin{align*}
\epsilon S &= S \epsilon \circ \delta\text, & \eta T &= \delta \circ T \eta\text, \\
S \Delta \circ \delta &= \delta T \circ T \delta \circ \Delta S\text, &
\text{and \ \ \quad} \delta \circ T \mu &= \mu T \circ S \delta \circ \delta
S\text.
\end{align*}
Then the two-sided Kleisli category $Kl(\delta)$ of the distributive law
$\delta$ has:
\begin{itemize*}
\item \textbf{Objects} those of $\C$;
\item \textbf{Maps} $A \to B$ in $Kl(\delta)$ given by maps $TA \to SB$ in
    $\C$,
\item \textbf{Identity maps} $\id_A \colon A \to A$ in $Kl(\delta)$ given by
    the map
\[TA \xrightarrow{\epsilon_A} A \xrightarrow{\eta_A} SA\]
in $\C$;
\item \textbf{Composition} for maps $f \colon A \to B$ and $g \colon B \to C$
    in $Kl(\delta)$ given by the map
\[TA \xrightarrow{\Delta_A} TTA \xrightarrow{Tf} TSB \xrightarrow{\delta_B} STB
\xrightarrow{Sg} SSC \xrightarrow{\mu_C} SC\] in $\C$.
\end{itemize*}
\end{Defn}
\noindent Now, we can emulate such a construction one dimension up:
\begin{Defn}
Let $\mathcal B$ be a bicategory, let $(S, \eta, \mu, \lambda, \rho, \tau)$ be
a pseudomonad and $(T, \epsilon, \Delta, \lambda', \rho', \tau')$ a
pseudocomonad on $\mathcal B$, and let $(\delta, \overline \eta, \overline
\epsilon, \overline \mu, \overline \Delta)$ be a pseudo-dis\-trib\-u\-tive law
of the pseudocomonad over the pseudomonad. Then the two-sided Kleisli
bicategory $Kl(\delta)$ of the pseudo-distributive law $\delta$ has:
\begin{itemize*}
\item \textbf{Objects} those of $\mathcal B$;
\item \textbf{Hom-categories} given by $Kl(\delta)(X, Y) = \mathcal B(TX, SY)$;
\item \textbf{Identity map} at $X$ given by the composite
\[ T X \xrightarrow{ \epsilon_X} X \xrightarrow{ \eta_X}  S X\text;\]
\item \textbf{Composition} $Kl(\delta)(Y, Z) \times Kl(\delta)(X, Y) \to
    Kl(\delta)(X, Z)$ given by
 \[\cd[@R+1em]{
  \mathcal B( T Y,  S Z) \times \mathcal B( T X,  S Y)
   \ar[d]^{\cong} \\
  1 \times \mathcal B( T Y,  S Z) \times 1 \times \mathcal B( T X,  S Y) \times 1
   \ar[d]^-{\elt{ \mu_Z} \times  S \times \elt{\delta_Y} \times  T \times \elt{ \epsilon_X}}\\
*++\txt<20pc>{$
  \mathcal B( S  S Z,  S Z) \times \mathcal B( S  T Y,
 S  S Z) \times \mathcal B( T  S Y,  S  T Y) \times \mathcal B( T
  T X,  T  S Y)
  \times \mathcal B( T X,  T  T X)$}
   \ar[d]^{\otimes}\\
  \mathcal B( T X,  S Z)}
 \]
where we use $\otimes$ to stand for some choice of order of composition for the
displayed fivefold composite. Explicitly, on maps, this composition is given by
taking for $(TY \xrightarrow{G} SZ) \otimes (T X \xrightarrow{F} SY) $ (some
choice of bracketing for) the composite
\[ TX
   \xrightarrow{\Delta_X}
   T  T X
   \xrightarrow{ T F}
   T  SY
   \xrightarrow{\delta_Y}
   S  T Y
   \xrightarrow{ S G}
   S SZ
   \xrightarrow{ \mu_Z}
   S Z\text.
\]
\end{itemize*}
\end{Defn}

Again, we shall not provide the associativity and unitality constraints
required to make this into a bicategory: they are now constructed from the
pseudomonad structure of $S$, the pseudocomonad structure of $T$ and the
pseudo-distributive structure of $\delta$.

Returning to the case under consideration, we claim that there is a
pseudo-distributive law $\delta \colon \hat S_c \hat S_m \Rightarrow \hat S_m
\hat S_c$ given as follows. Recall that we have $\hat S_c = \hat S_m = \hat S$,
and thus the component $\delta_{\cat C} \colon \hat S_c \hat S_m \cat C \tor
\hat S_m \hat S_c \cat C$ of $\delta$ is given by a functor $(SS\cat C)^\op
\times SS\cat C \to \cat{Set}$. So, given a \emph{discrete} category $X$, we
wish to take $\delta_{X}(\{\Sigma_m\}_{1 \leqslant m \leqslant j};
\{\Gamma_n\}_{1 \leqslant n \leqslant k})$ to be the set of suitable matchings
of $\{\Sigma_m\}$ with $\{\Gamma_n\}$. If we unwrap the definition of two-sided
Kleisli bicategory above, we now see that the desired monoidal structure on
$[(SX)^\op \times SX, \cat{Set}]$ is given precisely by horizontal composition
in $Kl(\delta)(X, X)$.

Thus we should \emph{like} to define a polycompositional polycategory to be a
monad on a discrete object $X$ in the bicategory $Kl(\delta)$; but to do this,
we must first establish the existence of the pseudo-distributive law $\delta$.
It is the task of the remainder of this paper to do this.

[The following alternative approach to the theory of polycategories was
suggested by Robin Houston: from the paper \cite{daystreet:lax},
multicategories with object set $X$ can be viewed as \emph{lax
  monoids} on the discrete object $X$ in $\cat{Mod}$. We might hope to
extend this to a notion of \emph{lax Frobenius algebra}, following
\cite{st:frob}; then a polycategory would be such a lax Frobenius algebra on a
discrete object of $\cat{Mod}$. However, we shall not pursue this further
here.]

\section{Deriving the pseudo-distributive law $\delta$}
We intend to construct the pseudo-distributive law $\delta$ by exploiting the
theory of double clubs, as developed in the companion paper \cite{rhgg1}. A
\emph{double club} is a generalisation of Kelly's abstract notion of club
\cite{kelly:clubs} from the level of categories to that of \emph{pseudo} (or
\emph{weak}) double categories. Let us recap briefly the details we shall need
here.

A \defn{pseudo double category} $\db K$ is a `pseudo-category' object in
$\cat{Cat}$. Explicitly, it consists of \defn{objects} $X, Y, Z, \dots$,
\defn{vertical maps} $f \colon X \to Y$, \defn{horizontal maps} $\b X \colon X_s \tor X_t$ and
\defn{cells}
\[\cd{
  X_s \ar[d]_{f_s} \ar[r]|-{\object@{|}}^{\X} \dtwocell{dr}{\b f}& X_t \ar[d]^{f_t}  \\
  Y_s \ar[r]|-{\object@{|}}_{\Y} & Y_t\text,}\]
together with notions of vertical and horizontal composition such that vertical
composition is associative on the nose, whilst horizontal composition is
associative up to invertible \emph{special cells}, where a cell is said to be
\emph{special} just when its vertical source and target maps are identities.
The objects and vertical maps of a pseudo double category $\db K$ form a
category $K_0$, whilst the horizontal maps and cells form a category $K_1$.

Any pseudo double category $\db K$ contains a bicategory $\mathcal B \db K$
consisting of the objects, horizontal maps and special cells of $\db K$; and it
is reasonable to think of $\db K$ as being the bicategory $\mathcal B\db K$
with `added vertical structure'. For example, we will be concerned with the
pseudo double category $\dCat$ which has:
\begin{itemize*}
\item \textbf{Objects} being small categories $\cat C$;
\item \textbf{Vertical maps} being functors $f \colon \cat C \to \cat D$;
\item \textbf{Horizontal maps} being profunctors $F \colon \cat D^\op \times \C
    \to
    \cat{Set}$; and
\item \textbf{Cells}
\[\cd{
\cat{C} \ar[d]_{f} \ar[r]|-{\object@{|}}^{H} \dtwocell{dr}{\alpha}& \cat{D} \ar[d]^{g}  \\
\cat{E} \ar[r]|-{\object@{|}}_{K} & \cat{F}}\] being natural transformations
\[\cd{
 \cat{D}^\op \times \cat{C}
  \ar[rr]^{g^\op \times f}
  \ar[dr]_{H} &
 {}
  \rtwocell{d}{\alpha} &
 \cat{F}^\op \times \cat{E}
  \ar[dl]^K \\ &
 \cat{Set}\text.
}\]
\end{itemize*}

Following the above philosophy, we think of $\dCat$ as being the bicategory
$\mathcal B \dCat = \cat{Mod}$ of categories, profunctors and profunctor
transformations, extended with the vertical structure of honest functors.

We can now go on to give a notion of \emph{homomorphism} of pseudo double
categories, extending that for bicategories, and two notions of
\emph{transformation} between homomorphisms, namely \emph{vertical} and
\emph{horizontal}: the former having vertical maps for its components, and the
latter horizontal. The correct notion of \emph{modification} for pseudo double
categories is that of a `cell' bordered by two horizontal and two vertical
transformations. In fact, it genuinely \emph{is} a cell in that we have:
\begin{Prop}
Given pseudo double categories $\db K$ and $\db L$, there is a pseudo double
category $[\db K, \db L]_\psi$ of homomorphisms $\db K \to \db L$, vertical
transformations, horizontal transformations and modifications.
\end{Prop}

 Pseudo double categories, homomorphisms and vertical transformations form themselves into a 2-category
$\cat{DblCat}_\psi$, and thus we can read off notions such as \emph{equivalence
of pseudo double categories} (equivalence in $\cat{DblCat}_\psi$) and
\emph{double monad} (monad in $\cat{DblCat}_\psi$).

We now recap very briefly the theory of double clubs developed in \cite{rhgg1}.
Given a homomorphism $S \colon \db K \to \db L$, we can construct the `slice
pseudo double category' $[\db K, \db L]_\psi / S\I$. It has
\begin{itemize*}
\item \textbf{Objects} $(A, \alpha)$ being homomorphisms $A \colon
  \db K \to \db L$ together with a vertical transformation $\alpha
  \colon A \Rightarrow S$;
\item \textbf{Vertical maps} $\gamma \colon (A, \alpha) \to (B, \beta)$
  being vertical transformations $\gamma \colon A \Rightarrow B$ such that $\beta\gamma = \alpha$;
\item \textbf{Horizontal maps} $(\b A, \b \alpha) \colon (A_s, \alpha_s) \tor
    (A_t, \alpha_t)$ being horizontal transformations $\b A \colon A_s \Tor
    A_t$ together with a modification
\[\cd{
A_s \ar@2[d]_{\alpha_s} \ar@2[r]|-{\object@{|}}^{\b A} \dthreecell{dr}{\b \alpha}& A_t \ar@2[d]^{\alpha_t}  \\
S \ar@2[r]|-{\object@{|}}_{S\I} & S\text.}\]
\item \textbf{Cells}
\[\cd{
  (A_s, \alpha_s) \ar[d]_{\gamma_s} \ar[r]|-{\object@{|}}^{(\b A, \b \alpha)} \dtwocell{dr}{\b \gamma}& (A_t, \alpha_t) \ar[d]^{\gamma_t}  \\
  (B_s, \beta_s) \ar[r]|-{\object@{|}}_{(\b B, \b \beta)} & (B_t, \beta_t)}\]
being modifications
\[\cd{
  A_s \ar@2[d]_{\gamma_s} \ar@2[r]|-{\object@{|}}^{\b A} \dthreecell{dr}{\b \gamma}& A_t \ar@2[d]^{\gamma_t}  \\
  B_s \ar@2[r]|-{\object@{|}}_{\b B} & B_t}\]
such that $\b \beta \b \gamma = \b \alpha$.
\end{itemize*}

For a sufficiently well-behaved $S$, this has a sub-pseudo double category
$\coll(S)$, whose objects are cartesian vertical transformations into $S$ and
whose horizontal maps are cartesian modifications into $S\I$. Here, a vertical
transformation or modification is said to be \emph{cartesian} just when all its
naturality squares are pullbacks; and so $\coll(S)$ is the pseudo double
category analogue of the `category of collections' $Coll(S)$ in Kelly's theory
of clubs.

We have a strict double homomorphism $\textsf{ev}_1 \colon \coll(S) \to \db L /
S\I_1$ which evaluates at $1$, where $1$ is the terminal object of $\db L$; and
as in the theory of clubs, we effectively lose no information in applying this
homomorphism:
\begin{Prop}\label{sufficiently}
For $\db L$ sufficiently complete, the strict double homomorphism
$\mathsf{ev}_1$ forms one side of an equivalence of pseudo double categories
\[\coll(S) \simeq \db L / S\I_1\text.\]
\end{Prop}
\begin{proof}
See \cite{rhgg1}, Proposition 30. \qedhere
\end{proof}

In order to give a sensible definition of `double club', we need a notion of
monoidal structure for pseudo double categories:
\begin{Defn}
A \defn{monoidal pseudo double category} is a pseudomonoid in the (cartesian)
monoidal 2-category $\cat{DblCat}_\psi$. \end{Defn}
\begin{Prop}
The `endohom' pseudo double category $[\db K, \db K]_\psi$ has a canonical
structure of monoidal pseudo double category; furthermore, given a double monad
$(S, \eta, \mu)$ on $\db K$, the slice pseudo double category $[\db K, \db
K]_\psi / S\I$ has a canonical monoidal structure lifting that of $[\db K, \db
K]_\psi$.
\end{Prop}
\begin{proof}
See \cite{rhgg1}, Propositions\ 39 \& 43. \qedhere
\end{proof}

We now have:
\begin{Defn}
A double monad $(S, \eta, \mu)$ on a pseudo double category $\db K$ is a
\defn{double club} if $\coll(S)$ is closed under the monoidal structure of
$[\db K, \db K]_\psi / S\I$.
\end{Defn}

Probably the best-known (and indeed, the oldest) example of a club is that for
\emph{symmetric strict monoidal categories} on $\cat{Cat}$. In \cite{rhgg1}, we
show that this club extends to a double club $(S, \eta, \mu)$ on $\dCat$; and
it is this result that we shall make use of in the rest of this section.

\subsection{Lifting to $\coll(S)$}
We wish to apply the theory of double clubs to simplifying the construction of
our pseudo-distributive law $\delta$.  Now, this pseudo-distributive law is
specified in terms of certain data and axioms in the bicategory $[\cat{Mod},
\cat{Mod}]_\psi$. However, it makes sense in \emph{any} bicategory equipped
with well-behaved notions of `whiskering' (\emph{well-behaved} in the sense
that they obey axioms formally similar to those for a $\cat{Gray}$-monoid
\cite{daystreet:monoidalbicats}).

We show in the Appendix of \cite{rhgg1} that for any double club, $\coll(S)$ is
not only a monoidal pseudo double category, but is equipped with a notion of
`whiskering', and it follows from this that $\mathcal B\big(\coll(S)\big)$ is a
suitable setting for the construction of a pseudo-distributive law.
Furthermore, it's easy see that there is a strict homomorphism of bicategories
\[V \colon \mathcal B\big(\coll(S)\big) \to \mathcal B\big([\dCat, \dCat]_\psi\big) \to [\cat{Mod}, \cat{Mod}]_\psi\]
which first forgets the projections onto $S\I$, and then forgets the vertical
structure; and moreover, that this homomorphism respects the `whiskering'
operations on these two bicategories. So if we can lift the pseudomonad $\hat
S_m$ and pseudocomonad $\hat S_c$ along $V$, then any pseudo-distributive law
we construct between their respective liftings will induce a
pseudo-distributive law between $\hat S_m$ and $\hat S_c$ as desired.

At this stage, it might appear that we have only made things more complicated,
by requiring ourselves to construct a pseudo-distributive law in $\coll(S)$;
but now we are in a position to utilise the equivalence of pseudo double
categories $\coll(S) \simeq \dCat / S\I_1$ in order to reduce the construction
of a pseudo-distributive law in $\coll(S)$ to a much simpler construction `at
$1$'.

So, let us begin by showing how we may lift our pseudomonad $\hat S_m$ and
pseudocomonad $\hat S_c$ to $\mathcal B \big(\coll(S)\big)$. The first stage is
straightforward; we lift
\[\hat S_m \ \ \text{to} \ \cd{S_m \ar@2[d]^{\id_S} \\ S}
 \qquad \text{and} \qquad \hat S_c \ \ \text{to} \ \cd{S_c \ar@2[d]^{\id_S} \\ S},\]
where, again, we are using $S_c$ and $S_m$ as aliases for $S \colon \dCat \to
\dCat$. Next we must lift $\hat \eta$, $\hat \mu$, $\hat \epsilon$ and $\hat
\Delta$ to horizontal transformations and cartesian modifications as follows:
\begin{gather*}
\cd[@+1em]{
 \id_{\dCat}
  \dthreecell{dr}{\tilde {\b \eta}}
  \ar@2[d]_{\eta}
  \ar@2[r]|-{\object@{|}}^-{\b \eta} &
 S_m
  \ar@2[d]^{\id_{S}} \\
 S
  \ar@2[r]|-{\object@{|}}_-{S\I} &
 S
}\ \text, \qquad \cd[@+1em]{
 S_mS_m
  \dthreecell{dr}{\tilde {\b \mu}}
  \ar@2[d]_{\mu}
  \ar@2[r]|-{\object@{|}}^-{\b \mu} &
 S_m
  \ar@2[d]^{\id_{S}} \\
 S
  \ar@2[r]|-{\object@{|}}_-{S\I} &
 S
}\ \text,\\ \cd[@+1em]{
 S_c
  \dthreecell{dr}{\tilde {\b \epsilon}}
  \ar@2[d]_{\id_{S}}
  \ar@2[r]|-{\object@{|}}^-{\b \epsilon} &
 \id_{\dCat}
  \ar@2[d]^{\eta} \\
 S
  \ar@2[r]|-{\object@{|}}_-{S\I} &
 S
}\qquad \text{and} \qquad \cd[@+1em]{
 S_c
  \dthreecell{dr}{\tilde {\b \Delta}}
  \ar@2[d]_{\id_{S}}
  \ar@2[r]|-{\object@{|}}^-{\b \Delta} &
 S_cS_c
  \ar@2[d]^{\mu} \\
 S
  \ar@2[r]|-{\object@{|}}_-{S\I} &
 S\text.
}
\end{gather*}

Now, to give the horizontal transformation $\b \eta$ we must give a `components
functor' $\dCat_0 \to \dCat_1$ along with `pseudonaturality' special cells. For
the former, we take the component at an object $X$ to be given by the component
of $\hat \eta$ at $X$, and the component at a vertical map $f$ to be given by
the pasting
\[\cd[@+2em]{
  X \dtwocell{dr}{\hat \eta_{(f_\ast)}^{-1}}
  \ar[r]|-{\object@{|}}^{\hat \eta_X} \ar[d]|-{\object@{|}}_{f_\ast} &
  SX \ar[d]|-{\object@{|}}^{\hat S(f_\ast)}
  \ar@/^3em/[d]|-{\object@{|}}^{(Sf)_\ast} \\
  Y \ar[r]|-{\object@{|}}_{\hat \eta_Y} & SY\text.  }\] For the
latter, we merely take the pseudonaturality 2-cells of $\hat \eta$; checking
all required naturality and coherence is now routine.  To give the cartesian
modification $\tilde {\b
  \eta}$, we must give components $\tilde {\b \eta}_X$ as follows:
\[\cd[@+1em]{
 X
  \dtwocell{dr}{\tilde {\b \eta}_X}
  \ar[d]_{\eta_X}
  \ar[r]|{\object@{|}}^{\hat \eta_X} &
 SX
  \ar[d]^{\id_{S}} \\
 SX
  \ar[r]|{\object@{|}}_{S\I_X} &
 SX\text.
}\] But this is to give natural families of maps $\hat \eta_X(y; x) \to
S\I_X(y; \spn{x})$ which we do via the natural isomorphisms \[\hat \eta_X(y; x)
\cong SX(y, \spn{x}) \cong S\I_X(y; \spn{x})\text.\] Checking naturality and
cartesianness is routine. We proceed similarly to lift $\hat \mu$, $\hat
\epsilon$ and $\hat \Delta$.

Finally, we must check that the modifications $\lambda$, $\rho$, $\tau$,
$\lambda'$, $\rho'$ and $\tau'$ for $\hat S_m$ and $\hat S_c$ lift to
$\coll(S)$. For example, we must check that
\[\lambda \colon \id_{\hat S} \Rrightarrow \hat \mu \otimes \hat S_m \hat \eta \colon \hat
S_m
\To \hat S_m\] lifts to a special modification
\[
 \b \lambda \colon \I_{(S_m, \id_S)}
  \Rrightarrow (\b \mu, \tilde{\b \mu}) \otimes (S_m, \id_S)(\b \eta, \tilde{\b \eta})
  \colon (S_m, \id_S) \Tor (S_m, \id_S)\text.
\]
This amounts to checking that the components of $\lambda$ are natural with
respect to cells of $\dCat$, and that they are compatible with the projections
down to $S\I$; and this is merely a matter of diagram chasing.

Therefore, in order to obtain our desired pseudo-distributive law on
$\cat{Mod}$, it suffices to produce data and axioms for a pseudo-distributive
law between $(S_m, \id_S)$ and $(S_c, \id_S)$ as detailed above. We now wish to
see how we can use the theory of double clubs to reduce this to data and axioms
in $\dCat / S\I_1$.

\subsection{Reducing to $\dCat / S\I_1$}
We begin with (PDD1), for which we must produce a horizontal arrow
\[(\b \delta, \tilde{\b \delta}) \colon (S_cS_m, \mu) \tor (S_mS_c, \mu)\]
of $\coll(S)$, i.e., a horizontal transformation and a cartesian modification
as follows:
\[\cd{
 S_cS_m
  \dthreecell{dr}{\tilde{\b \delta}}
  \ar@2[d]_{\mu}
  \ar@2[r]|-{\object@{|}}^{\b \delta} &
 S_mS_c
  \ar@2[d]^{\mu} \\
 S
  \ar@2[r]|-{\object@{|}}_{S\I} &
 S\text.
}\] Now, suppose we have a horizontal arrow
\[\cd{
 S_cS_m1
  \dtwocell{dr}{\tilde{\b \delta}_1}
  \ar[d]_{\mu_1}
  \ar[r]|-{\object@{|}}^{{\b \delta}_1} &
 S_mS_c1
  \ar[d]^{\mu_1} \\
 S1
  \ar[r]|-{\object@{|}}_{S\I_1} &
 S1
}\] of $\dCat / S\I_1$. We should like to say that $({\b \delta}_1, \tilde{\b
\delta}_1)$ is the component at $1$ of some horizontal arrow $(\b \delta,
\tilde{\b \delta})$ of $\coll(S)$, which amounts to asking for the double
homomorphism $\mathsf{ev}_1 \colon \coll(S) \to \dCat / S\I_1$ to be
`horizontally full', in the following sense:
\begin{Prop}\label{horzfull}
Let $(A_s, \alpha_s)$ and $(A_t, \alpha_t)$ be objects of $\coll(S)$, and
suppose that we have a horizontal arrow
\[\cd{
 A_s1
  \dtwocell{dr}{\b \gamma}
  \ar[d]_{(\alpha_s)_1}
  \ar[r]|-{\object@{|}}^{\b a} &
 A_t1
  \ar[d]^{(\alpha_t)_1} \\
 S1
  \ar[r]|-{\object@{|}}_{S\I_1} &
 S1
}\] of $\dCat / S\I_1$. Then there is a horizontal arrow $(\b A, \b \Gamma)$ of
$\coll(S)$:-
\[\cd{
 A_s
  \dthreecell{dr}{\b \Gamma}
  \ar@2[d]_{\alpha_s}
  \ar@2[r]|-{\object@{|}}^{\b A} &
 A_t
  \ar@2[d]^{\alpha_t} \\
 S
  \ar@2[r]|-{\object@{|}}_{S\I} &
 S
}\] such that $\mathsf{ev}_1(\b A, \b \Gamma) = (\b a, \b \gamma)$.
\end{Prop}
\begin{proof}
Proposition \ref{sufficiently} above tells us that $\mathsf{ev}_1 \colon
\coll(S) \to \dCat / S\I_1$ forms one side of an equivalence of double
categories: and the proof of this given in \cite{rhgg1} constructs an explicit
quasi-inverse $\dCat / S\I_1 \to \coll(S)$. The following is a simple
adaptation of this construction to the problem at hand.

To give the horizontal transformation $\b A$, we must give, amongst other
things, a component profunctor $\b A_\C \colon A_s\C \tor A_t\C$ at each small
category $\C$; whilst to give $\b \Gamma$, we must give, for each small
category $\b C$, a cell
\[\cd{
 A_s\C
  \dtwocell{dr}{\b \Gamma_\C}
  \ar[d]_{(\alpha_s)_\C}
  \ar[r]|-{\object@{|}}^{\b A_\C} &
 A_t\C
  \ar[d]^{(\alpha_t)_\C} \\
 S\C
  \ar[r]|-{\object@{|}}_{S\I_\C} &
 S\C
}\] of $\dCat$. We may view $\b \gamma$ as a morphism $\b a \to S\I_1$ in
the category $\dCat_1$ of profunctors and transformations between them; and
thus may form $\b A_\C$ and $\b \Gamma_\C$ as the following pullback in
$\dCat_1$:
\begin{equation}\label{specialsquare}\tag{*}
    \cd{
        \b A_\C \ar[r] \ar[d]_{\b \Gamma_\C} & \b a \ar[d]^{\b \gamma} \\
        S\I_\C \ar[r]_{S\I_!} & S\I_1\text.
    }
\end{equation}

Now, in order that $\b A_\C$ and $\b \Gamma_\C$ should have the correct sources
and targets, we must choose the pullback (*) in such a way that application of
the source and target functors $s, t \colon \dCat_1 \to \dCat_0$ sends it to
the respective squares:
\begin{equation*}
    \cd{
        A_s\C \ar[r] \ar[d]_{(\alpha_s)_\C} & A_s1 \ar[d]^{(\alpha_s)_1} \\
        S\C \ar[r]_{S!} & S1
    } \qquad \text{and} \qquad
    \cd{
        A_t\C \ar[r] \ar[d]_{(\alpha_t)_\C} & A_t1 \ar[d]^{(\alpha_t)_1} \\
        S\C \ar[r]_{S!} & S1
    }
\end{equation*}
in $\dCat_0$. That we may do this follows from two observations: firstly, that
both the displayed squares are pullbacks in $\dCat_0$, by cartesianness of
$\alpha_s$ and $\alpha_t$; and secondly, that the functor $(s, t) \colon
\dCat_1 \to \dCat_0 \times \dCat_0$ (strictly) creates pullbacks.

In order that we should have $\textsf{ev}_1(\b A, \b \Gamma) = (\b a, \b
\gamma)$, we make one further demand: that when $\C = 1$, the pullback square
(*) should be chosen as
\begin{equation*}
    \cd{
        \b a \ar[r]^{\id} \ar[d]_{\b \gamma} & \b a \ar[d]^{\b \gamma} \\
        S\I_1 \ar[r]_{\id} & S\I_1\text.
    }
\end{equation*}

Apart from this care in choosing the pullback squares (*), the remaining
details in the construction of $\b A$ and $\b \Gamma$ are exactly as in the
proof of Proposition 30 of \cite{rhgg1}, and hence omitted. \qedhere
\end{proof}

Thus, given a horizontal arrow $({\b \delta}_1, \tilde{\b \delta}_1) \colon
(S_cS_m1, \mu_1) \tor (S_mS_c1, \mu_1)$ of $\dCat / S\I_1$, we can produce a
horizontal arrow $(\b \delta, \tilde{\b \delta}) \colon (S_cS_m, \mu) \tor
(S_mS_c, \mu)$ of $\coll(S)$ whose image under $\mathsf{ev}_1$ is precisely
$({\b \delta}_1, \tilde{\b \delta}_1)$.

To derive the remaining data (PDD2) and (PDD3), we observe the following: the
double homomorphism $F := \mathsf{ev}_1 \colon \coll(S) \to \dCat / S\I_1$ is
built upon two functors $F_0 \colon \coll(S)_0 \to \big(\dCat / S\I_1\big)_0$
and $F_1 \colon \coll(S)_1 \to \big(\dCat / S\I_1\big)_1$; and since $F$ forms
one side of an equivalence of pseudo double categories, it follows that $F_0$
and $F_1$ each form one side of an equivalence of ordinary categories.  In
particular, the functor $F_1 \colon \coll(S)_1 \to \dCat_1 / S\I_1$ is full and
faithful. Thus, considering $\overline {\b \eta}$ for instance, we must find a
special invertible cell
\[\overline {\b \eta} \colon (\b \delta, \tilde{\b \delta})
\otimes (S_m, \id_S)(\b \eta, \tilde{\b \eta}) \Rrightarrow (\b
\eta, \tilde{\b \eta}) (S_m, \id_S)\] of $\coll(S)$. Since $F_1$ is full
and faithful, it suffices for this to find a special invertible cell
\[\overline {\b \eta}_1 \colon ({\b \delta}_1, \tilde{\b \delta}_1)
\otimes \big((S_m, \id_{S})(\b \eta, \tilde{\b \eta})\big)_1
\Rightarrow \big((\b \eta, \tilde{\b \eta}) (S_m, \id_{S})\big)_1\]
of $\dCat / S\I_1$. We proceed similarly for the remaining data.

Finally, we must ensure that (PDA1)--(PDA10) are satisfied, which amounts to
checking certain equalities of pastings in $\mathcal B \big(\coll(S)\big)$,
which in turn amounts to checking certain equalities of maps in $\coll(S)_1$;
but since the functor $F_1 \colon \coll(S)_1 \to \dCat_1 / S\I_1$ is faithful,
it suffices for this to check that these equalities hold in $\dCat / S\I_1$.

\section{Constructing the pseudo-distributive law at $1$}
\subsection{The double club $S$ on $\dCat$}\label{presentation}
In order to construct the data and axioms laid out at the end of the previous
section, we will require a detailed presentation of the double club $(S, \eta,
\mu)$ on $\dCat$. Since this double club looks like the free symmetric monoidal
category monad on $\cat{Cat}$ in the vertical direction, and like its lifting
$\hat S$ to $\cat{Mod}$ in the horizontal direction, we may do this by giving a
presentation of these latter two entities.
\begin{Defn}
We write $S1$ for the category of finite cardinals and bijections, with:
\begin{itemize*}
\item \textbf{Objects} the natural numbers $0, 1, 2, \dots$;
\item \textbf{Maps} $\sigma \colon n \to m$ bijections of $\{1, \dots, n\}$
    with $\{1, \dots, m\}$,
\end{itemize*}
and with composition and identities given in the evident way.
\end{Defn}

\begin{Defn}
The \emph{free symmetric strict monoidal category} 2-functor $S \colon
\cat{Cat} \to \cat{Cat}$ is given as follows:
\begin{itemize*}
\item \textbf{On objects}: Given a small category $\cat C$, we give $S\cat C$
    as follows:
\begin{itemize*}
\item \textbf{Objects} of $S\cat C$ are pairs $(n, \spn{c_i})$, where $n \in
    S1$ and $c_1, \dots, c_n \in \ob \cat C$;
\item \textbf{Arrows} of $S\cat C$ are
\[(\sigma, \spn{g_i}) \colon (n, \spn{c_i}) \to (m, \spn{d_i})\text,\] where
$\sigma \in S1(n, m)$ and $g_i \colon c_i \to d_{\sigma(i)}$ (note that
necessarily $n = m$).
\end{itemize*}
Composition and identities in $S\cat C$ are given in the evident way; namely,
\begin{align*}
\id_{(n, \spn{c_i})} &= (\id_n, \spn{\id_{c_i}})\\
\text{and } (\tau, \spn{g_i}) \circ (\sigma, \spn{f_i}) &= (\tau
\circ \sigma, \spn{g_{\sigma(i)} \circ f_i})\text.
\end{align*}
\item \textbf{On maps}: Given a functor $F \colon \cat C \to \cat D$, we give
    $SF \colon S\cat C \to S\cat D$ by
\[
SF(n, \spn{c_i}) = (n, \spn{Fc_i})\quad \text{ and } \quad
SF(\sigma, \spn{g_i}) = (\sigma, \spn{Fg_i})\text.
\]
\item \textbf{On 2-cells}: Given a natural transformation $\alpha \colon F
    \Rightarrow G \colon \cat C \to \cat D$, we give $S\alpha \colon SF
    \Rightarrow SG \colon S\cat C \to S\cat D$ by
\[(S\alpha)_{(n, \spn{c_i})} = (\id_{n}, \spn{\alpha_{c_i}})\text.\]
\end{itemize*}
\end{Defn}
\noindent Now, although the above is sufficient to describe the iterated
functor $S^2 \colon \cat{Cat} \to \cat{Cat}$, it will be much more pleasant to
work with the following alternative presentation. First note that we may
describe $S^21$ as follows:

\begin{itemize}
\item \textbf{Objects} are order-preserving maps $\phi \colon n_\phi \to
    m_\phi$, where $n_\phi$, $m_\phi \in \mathbb N$. We write such an
    object simply as $\phi$, with the convention that $\phi$ has domain and
    codomain $n_\phi$ and $m_\phi$ respectively.
\item \textbf{Maps} $f \colon \phi \to \psi$ are pairs of bijections $f_n
    \colon n_\phi \to n_\psi$ and $f_m \colon m_\phi \to m_\psi$ such that
    the following diagram commutes:
\[
\cd[@!@-1em]{
  n_\phi \ar[d]_{\phi} \ar[r]^{f_n} & n_\psi \ar[d]^{\psi} \\
  m_\phi \ar[r]_{f_m} & m_\psi\text.
}
\]
\end{itemize}
\noindent It may not be immediately obvious that this \emph{is} a presentation
of $S^21$. The picture is as follows: an object $\phi$ of $S^21$ is to be
thought of as a collection of $n_\phi$ points partitioned into $m_\phi$ parts
in accordance with $\phi$. Given such an object, one can permute internally any
of its $m_\phi$ parts, or can in fact permute the set of $m_\phi$ parts itself;
and a typical map describes such a permutation. For example, the objects
\begin{align*}\phi \colon 5 &\to 4 & \psi \colon 5 &\to 4\\
1, 2, 3, 4, 5 &\mapsto 1, 1, 3, 4, 4 & 1, 2, 3, 4, 5 &\mapsto 2,
2, 3, 4, 4
\end{align*}
should be visualised as
\[
\dbr{\dbr{\bullet, \bullet}, \dbr{}, \dbr{\bullet}, \dbr{\bullet,
    \bullet}} \quad \text{\ \ and\ \ } \quad \dbr{\dbr{},
  \dbr{\bullet, \bullet}, \dbr{\bullet}, \dbr{\bullet, \bullet}}
\]
respectively, whilst a typical map $\phi \to \psi$ is given by
\begin{align*}f_n \colon 5 &\to 5 & f_m \colon 4 &\to 4\\
1, 2, 3, 4, 5 &\mapsto 5, 4, 3, 1, 2 & 1, 2, 3, 4 &\mapsto 4, 1,
3, 2
\end{align*}
and should be visualised as
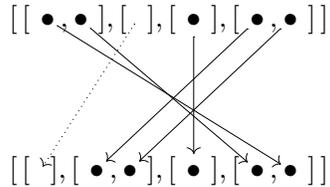
\[
\begin{xy}
\xymatrix"*"@1@C=0pt{
  [\,[\,&\bullet &, &\bullet&\,], [\,&&\,], [\,&\bullet&\,], [\,&\bullet&, &\bullet&\,]\,]
}\POS-(0,20) \xymatrix@1@C=0pt{
  [\,[\,&\ar@{<.}["*"rrrr]&\,], [\,&\bullet \ar@{<-}["*"rrrrrr]&, &\bullet \ar@{<-}["*"rrrrrr]&\,], [\,&\bullet \ar@{<-}["*"]&\,], [\,&\bullet\ar@{<-}["*"llllll]&, &\bullet\ar@{<-}["*"llllllllll]&\,]\,]\text.
}
\end{xy}
\]
So now, given a category $\cat C$, we can present $S^2\cat C$ as follows:
\begin{itemize}
\item \textbf{Objects} of $S^2\cat C$ are pairs $(\phi, \spn{c_i})$, where
    $\phi = n_\phi \to m_\phi \in S^21$ and $c_1, \dots, c_{n_\phi} \in \ob
    \cat C$;
\item \textbf{Arrows} of $S^2\cat C$ are
\[(f, \spn{g_i}) \colon (\phi, \spn{c_i}) \to (\psi, \spn{d_i})\text,\] where
$f = (f_n, f_m) \in S^21(\phi, \psi)$ and $g_i \colon c_i \to d_{f_n(i)}$;
composition and identities are given analogously to before.
\end{itemize}
\noindent We can extend the above in the obvious way to 1- and 2-cells of
$\cat{Cat}$ to give a presentation of the 2-functor $S^2$. Using this alternate
presentation of $S^2$, we may describe the rest of the 2-monad structure of
$S$:
\begin{Defn}
The 2-natural transformation $\eta \colon \id_{\cat{Cat}} \Rightarrow S$ has
component at $\cat C$ given by
\begin{align*}
\eta_{\cat C} \colon \cat C & \to S\cat C\\
x & \mapsto (1, \spn{x})\\
f & \mapsto (\id_1, \spn{f})\text,
\end{align*}
whilst the 2-natural transformation $\mu \colon S^2 \Rightarrow S$ has
component at $\cat C$ given by
\begin{align*}
\eta_{\cat C} \colon SS\cat C & \to S\cat C\\
(\phi, \spn{c_i}) & \mapsto (n_\phi, \spn{c_i})\\
(f, \spn{g_i}) & \mapsto (f_n, \spn{g_i})\text.
\end{align*}
\end{Defn}
\noindent We will also need to make use of the threefold iterate $S^3$, and so
it will be useful to present it in the above style. We first give $S^31$ as
follows:
\begin{itemize}
\item \textbf{Objects} are diagrams $\phi = n_\phi \xrightarrow{\phi_1}
    m_\phi \xrightarrow{\phi_2} r_\phi$ in the category of finite ordinals
    and order preserving maps;
\item \textbf{Maps} $f \colon \phi \to \psi$ are triples $(f_n, f_m, f_r)$
    of bijections making
\[
\cd[@!@-1em]{
  n_\phi \ar[d]_{\phi_1} \ar[r]^{f_n} & n_\psi \ar[d]^{\psi_1}\\
  m_\phi \ar[d]_{\phi_2} \ar[r]^{f_m} & m_\psi \ar[d]^{\psi_2} \\
  r_\phi \ar[r]_{f_r} & r_\psi \text.
}
\]
commute.
\end{itemize}
Whereupon we may describe $S^3\cat C$ as follows:
\begin{itemize}
\item \textbf{Objects} are pairs $(\phi, \spn{c_i})$, where $\phi = n_\phi
    \to m_\phi \to r_\phi \in S^31$ and $c_1, \dots, c_{n_\phi} \in \ob
    \cat C$;
\item \textbf{Arrows} are
\[(f, \spn{g_i}) \colon (\phi, \spn{c_i}) \to (\psi, \spn{d_i})\text,\] where
$f = (f_n, f_m, f_r) \in S^31(\phi, \psi)$ and $g_i \colon c_i \to
d_{f_n(i)}$.
\end{itemize}
\noindent As before, we may straightforwardly extend this definition to 1- and
2-cells of $\cat{Cat}$.

Finally, we give a presentation of the pseudomonad $(\hat S, \hat \eta, \hat
\mu)$ on $\cat{Mod}$:
\begin{Defn}
The homomorphism $\hat S \colon \cat{Mod} \to \cat{Mod}$ is given as follows:
\begin{itemize*}
\item \textbf{On objects}: Given a small category $\cat C$, we take $\hat S\cat
    C = S\cat C$;
\item \textbf{On maps}: Given a map $F \colon \cat C \tor \cat D$, the map
    $\hat SF \colon S\cat C \tor S\cat D$ is the following profunctor: an
    element of $\hat SF\big((n, \spn{d_i}); (m, \spn{c_i})\big)$ is given by
\[(\sigma, \spn{g_i}) \colon (n, \spn{d_i}) \tor (m, \spn{c_i})\text,\] where
$\sigma \in S1(n, m)$ and $g_i \in F(d_i; c_{\sigma(i)})$, whilst the action on
these elements by maps $(\tau, \spn{h_i}) \colon (m, \spn{c_i}) \to (m',
\spn{c'_i})$ and $(\upsilon, \spn{f_i}) \colon (n', \spn{d'_i}) \to (n,
\spn{d_i})$ is given by
\begin{align*}
(\sigma, \spn{g_i}) \cdot (\upsilon, \spn{f_i}) &= (\sigma \circ
\upsilon,
\langle g_{\upsilon(i)} \cdot f_i\rangle)\\
(\tau, \spn{h_i}) \cdot (\sigma, \spn{g_i}) &= (\tau \circ
\sigma, \langle h_{\sigma(i)} \cdot g_i\rangle )\text;
\end{align*}
\item \textbf{On 2-cells}: Given a transformation $\alpha \colon F \Rightarrow
    G \colon \cat C \tor \cat D$, we give $S\alpha \colon SF \Rightarrow SG
    \colon S\cat C \tor S\cat D$ by
\[(S\alpha)(\sigma, \spn{g_i}) = (\sigma, \spn{\alpha(g_i)})\text.\]
\end{itemize*}
Further, the pseudo-natural transformations
\begin{align*}
\hat \eta &\colon \id \Rightarrow \hat S \colon
\cat{Mod} \to \cat{Mod}\\
\text{and } \hat \mu &\colon \hat S^2 \Rightarrow \hat S \colon
\cat{Mod} \to \cat{Mod}
\end{align*}
have respective components
\[\hat \eta_X = (\eta_X)_\ast \quad \text{and} \quad \hat \mu_X =
(\mu_X)_\ast\text.\]
\end{Defn}

\subsection{Spans}
We shall also need a few preliminaries about acyclic and connected graphs. We
seek to capture their combinatorial essence in a categorical manner, allowing a
smooth presentation of the somewhat involved proof which follows.
%

The objects of our attention are spans in $\cat{FinCard}$, i.e., diagrams $n
\leftarrow k \rightarrow m$ in the category of finite cardinals and all maps.
When we write `span' in future, it should be read as `span in $\cat{FinCard}$'
unless otherwise stated. We also make use without comment of the evident
inclusions $\cat{FinOrd} \rightarrow \cat{FinCard}$ and $S1 \rightarrow
\cat{FinCard}$.

Now, each span $n \leftarrow k \rightarrow m$ determines a (categorist's) graph
$k \rightrightarrows n + m$; if we forget the orientation of the edges of this
graph, we get a (combinatorialist's) undirected multigraph. We say that a span
$n \leftarrow k \rightarrow m$ is
\defn{acyclic} or \defn{connected} if the associated multigraph is so.
Note that the \emph{acyclic} condition includes the assertion that there are no
multiple edges.

\begin{Prop}
Given a span $n \xleftarrow{\theta_1} k \xrightarrow{\theta_2} m$, the number
of connected components of the graph induced by the span is given by the
cardinality of $r$ in the pushout diagram
\[\cd{
  k \ar[r]^{\theta_2} \ar[d]_{\theta_1} & m \ar[d]^{\tau_2} \\
  n \ar[r]_{\tau_1} & r
}\] in $\cat{FinCard}$.
\end{Prop}
\begin{proof}
  Given the above pushout diagram, set $n_i = \tau_1^{-1}(i)$ and $m_i
  = \tau_2^{-1}(i)$ (for $i = 1, \dots, r$).  Now we observe that, for
  $i \neq j$, we have
\[\theta_1^{-1}(n_i) \cap \theta_2^{-1}(m_j) = \theta_1^{-1}(n_i) \cap \theta_1^{-1}(n_j) = \emptyset\text,\]
so that induced graph of the span has at least $r$ unconnected parts (with
respective vertex sets $n_i + m_i$). On the other hand, if the induced graph
$G$ had strictly more than $r$ connected components, we could find vertex sets
$v_1, \dots, v_{r+1}$ which partition $v(G)$, and for which
\[x \in v_i\text,\  y \in v_j\text{ (for }i \neq j\text) \qquad \text{implies}
\qquad x \text{ is not adjacent to } y\text.
\eqno{(\text{\dag})}\] But now define maps $\tau_1 \colon n \to
r+1$ and $\tau_2 \colon m \to r+1$ by letting $\tau_i(x)$ be the $p$ for which
$x \in v_p$. Then by condition $(\text{\dag})$, we have $\tau_1(\theta_1(a)) =
\tau_2(\theta_2(a))$ for all $a \in k$, and so we have a commuting diagram
\[\cd{
  k \ar[r]^{\theta_2} \ar[d]_{\theta_1} & m \ar[d]^{\tau_2} \\
  n \ar[r]_-{\tau_1} & r+1
}\] for which the bottom right vertex does not factor through $r$,
contradicting the assumption that $r$ was a pushout. Hence $G$ has precisely
$r$ connected components. \qedhere
\end{proof}

\begin{Cor}
A  span $n \xleftarrow{\theta_1} k \xrightarrow{\theta_2} m$ is connected if
and only if the diagram
\[\cd{
  k \ar[r]^{\theta_2} \ar[d]_{\theta_1} & m \ar[d] \\
  n \ar[r] & 1
}\] is a pushout in $\cat{FinCard}$.
\end{Cor}

\begin{Prop}
A span $n \xleftarrow{\theta_1} k \xrightarrow{\theta_2} m$ is acyclic if and
only, for every monomorphism $\iota \colon k' \hookrightarrow k$,
\[
\cd{
  k \ar[r]^{\theta_2} \ar[d]_{\theta_1} & m \ar[d] \\
  n \ar[r] & r
} \quad \text{ a pushout implies } \quad \cd{
  k' \ar[r]^{\theta_2 \iota} \ar[d]_{\theta_1 \iota} & m \ar[d] \\
  n \ar[r] & r
}\text{ not a pushout.}\]
\end{Prop}
\begin{proof}
  Suppose the left hand diagram is a pushout; then the associated
  graph $G$ of the span has $r$ connected components.

  Suppose first that $G$ is acyclic, and $\iota \colon k' \hookrightarrow
  k$. Then the graph $G'$ associated to the span $n \xleftarrow{\theta_1 \iota}
  k' \xrightarrow{\theta_2 \iota} m$ has the same vertices as $G$ but
  strictly fewer edges; and since $G$ is acyclic, $G'$ must have strictly more than $r$
  connected components, and hence $r$ cannot be a pushout for the
  right-hand diagram.

  Conversely, if $G$ has a cycle, then we can remove some edge of
  $G$ without changing the number of connected components; and thus we
  obtain some monomorphism $\iota \colon k' \hookrightarrow k$ making
  the right-hand diagram a pushout. \qedhere
\end{proof}

\begin{Prop}
Suppose we have a commuting diagram
\[\cd{
  k \ar[r]^{\theta_2} \ar[d]_{\theta_1} & m \ar[d]^{\phi_2} \\
  n \ar[r]_{\phi_1} & r\text.
}\eqno{(*)}\] Then the spans $m^{(i)} \leftarrow k^{(i)}
\rightarrow n^{(i)}$ (for $i = 1, \dots, r$) induced by pulling back along
elements $i \colon 1 \to r$ are all connected if and only if $(*)$ is a
pushout.
\end{Prop}
\begin{proof}
Suppose all the induced spans are connected; then each diagram
\[\cd{
  k^{(i)} \ar[r]^{\theta_2^{(i)}} \ar[d]_{\theta_1^{(i)}} & m^{(i)} \ar[d] \\
  n^{(i)} \ar[r] & 1
}\] is a pushout; hence the diagram
\[\cd{
  {\sum_i} k^{(i)} \ar[r]^{{\sum_i} \theta_2^{(i)}} \ar[d]_{{\sum_i} \theta_1^{(i)}} & {\sum_i} m^{(i)} \ar[d] \\
  {\sum_i} n^{(i)} \ar[r] & r
}\] is also a pushout, whence it follows that $(*)$ is itself a
pushout.

Conversely, if $(*)$ is a pushout, then pulling this back along the map $i
\colon 1 \to r$ yields another pushout in $\cat{FinCard}$, so that each induced
span is connected. \qedhere
\end{proof}

\begin{Prop}
Let $G$ be a graph with finite edge and vertex sets. Any two of the following
conditions implies the third:
\begin{itemize}
\item $G$ is acyclic;
\item $G$ is connected;
\item $\abs{v(G)} = \abs{e(G)} + 1$.
\end{itemize}
\end{Prop}
\begin{proof}\hfill
\begin{itemize}
\item If $G$ is acyclic and connected, then it is a tree, and so
    $\abs{v(G)} = \abs{e(G)} + 1$;
\item if $G$ is connected with $\abs{v(G)} = \abs{e(G)} + 1$, then it is
    minimally connected, hence a tree, and so acyclic;
\item if $G$ is acyclic with $\abs{v(G)} = \abs{e(G)} + 1$, then it is
    maximally acyclic, hence a tree, and so connected. \qedhere
\end{itemize}
\end{proof}

\begin{Cor}
A span $n \xleftarrow{\theta_1} k \xrightarrow{\theta_2} m$ is acyclic and
connected if and only if the diagram
\[\cd{
  k \ar[r]^{\theta_2} \ar[d]_{\theta_1} & m \ar[d] \\
  n \ar[r] & 1
}\] is a pushout in $\cat{FinCard}$, and $n + m = k + 1$.
\end{Cor}

\begin{Cor}
Let there be given a commuting diagram
\[\cd{
 k \ar[r]^{\theta_2} \ar[d]_{\theta_1} & m \ar[d]^{\phi_2} \\
  n \ar[r]_{\phi_1} & r\text;
}\eqno{(*)}\] then the induced spans $m^{(i)} \leftarrow k^{(i)}
\rightarrow n^{(i)}$ (for $i = 1, \dots, r$) are acyclic and connected if and
only if $(*)$ is a pushout and $m + n = k + r$.
\end{Cor}

\subsection{(PDD1)}\label{pdd1sec}
We are now ready to give our pseudo-distributive law at $1$, and we begin with
(PDD1), for which we must give a horizontal arrow
\[\cd{
 S_cS_m1
  \dtwocell{dr}{\tilde{\b \delta}_1}
  \ar[d]_{\mu_1}
  \ar[r]|-{\object@{|}}^{\b \delta_1} &
 S_mS_c1
  \ar[d]^{\mu_1} \\
 S1
  \ar[r]|-{\object@{|}}_{S\I_1} &
 S1
}\] of $\dCat / S\I_1$.
\begin{Defn}\label{suitable}
The \defn{profunctor of suitable matchings}, $\b \delta_1 \colon \hat S_c \hat
S_m 1 \tor \hat S_m \hat S_c 1$ is the following functor $(S^21)^\op \times
S^21 \to \cat{Set}$:
\begin{itemize}
\item \textbf{On objects}: elements $f \in \b \delta_1(\phi; \psi)$ are
    bijections $f_n$ fitting into the diagram
\[
\cd[@!@-1em]{
  n_\phi \ar[r]^{f_n} \ar[d]_{\phi} & n_\psi \ar[d]^\psi \\
  m_\phi  & m_\psi
}
\]
such that the span $m_\phi \xleftarrow{\phi} n_\phi \xrightarrow{\psi \circ
f_n} m_\psi$ is acyclic and connected.

\item \textbf{On maps}: Let $g \colon \psi \to \rho$ in $S^21$ and let $f
    \in \b \delta_1(\phi; \psi)$. Then we give $g \cdot f \in \b
    \delta_1(\phi; \rho)$ by
\[
\cd[@!]{
  n_\phi \ar[r]^{g_n \cdot f_n} \ar[d]_{\phi} & n_\rho \ar[d]^{\rho} \\
  m_\phi & m_\rho
}
\]
This action is evidently functorial, but we still need to check that it
really does yield an element of $\b \delta_1(\phi; \rho)$; that is, we need
the associated span to be acyclic and connected. But this span is the top
path of the diagram
\[\cd[@!@-=2em]{& & n_\phi \ar[dr]^{f_n} \ar[ddll]_{\phi} \\
      & & & n_\psi \ar[dr]_{\psi} \ar[r]^{g_n} & n_\rho \ar[dr]^{\rho} \\
      m_\phi & & & & m_\psi \ar[r]_{g_m} & m_\rho\text;
}
\]
and therefore also the bottom path, since the right-hand square commutes.
But since $g_m$ is an isomorphism, the graph induced by the span $m_\phi
\xleftarrow{\phi} n_\phi \xrightarrow{\psi f_n} m_\psi$ is isomorphic to
the graph induced by the span $m_\phi \xleftarrow{\phi} n_\phi
\xrightarrow{g_m\psi f_n} m_\rho$, and hence the latter is acyclic and
connected since the former is. So we have a well-defined left action of
$S^21$ on $\b \delta_1$; and we proceed similarly to define an action on
the right.
\end{itemize}
\end{Defn}
We now give the 2-cell $ \tilde {\b \delta}_1$, for which we must give natural
families of maps $\b \delta_1(\phi; \psi) \to S1(n_\phi, n_\psi)$. But this is
straightforward: we simply send
\[
\cd[@!@-1em]{
  n_\phi \ar[r]^{f_n} \ar[d]_{\phi} & n_\psi \ar[d]^\psi \\
  m_\phi  & m_\psi
}
\]
in $\b \delta_1(\phi; \psi)$ to $f_n$ in $S1(n_\phi; n_\psi)$. It is visibly
the case that this satisfies the required naturality conditions.

Now, consider the transformation $\delta \colon \hat S_c \hat S_m \Rightarrow
\hat S_m \hat S_c$ induced by this $(\b \delta_1, \tilde{\b \delta}_1)$. From
Definition \ref{suitable} and Proposition \ref{horzfull}, we obtain that the
component of $\delta$ at a discrete category $X$ is given by
\begin{equation*}
    \delta_{X}(\{\Sigma_m\}; \{\Gamma_n\}) = \set{\sigma}{\text{$\sigma$ is a suitable matching of $\{\Sigma_m\}$ with $\{\Gamma_n\}$}}
\end{equation*}
as desired.


\subsection{(PDD2)}
For (PDD2) we must produce the component of the invertible special
modifications $\overline {\b \eta}$ and $\overline {\b \epsilon}$ at $1$:
\begin{Prop}
There is an invertible special cell
\[\cd[@+1.3em]{
 S_c1
  \ar[d]|-{\object@{|}}_{(S_c\b \eta)_1}
  \ar[dr]|-{\object@{|}}^{(\b \eta S_c)_1}
  & {} \\
 S_cS_m1
  \rtwocell[0.35]{ur}{\overline {\b \eta}_1}
  \ar[r]|-{\object@{|}}_{\b \delta_1} &
 S_mS_c1
}\] mediating the centre of this diagram in $\coll(S)$ (where we
omit the projections to $S\I$).
\end{Prop}

\begin{proof}
With respect to the descriptions of $S1$ and $S^21$ given above, we observe
that that the functors $(S\eta)_1 \colon S1 \to S^21$ and $\eta_{S1} \colon S1
\to S^21$ are given by
\begin{align*}
S \eta_1 \colon n & \mapsto (n \xrightarrow{\id} n) & \eta_{S1} \colon n & \mapsto (n \xrightarrow{!} 1) \\
f & \mapsto (f, f) & f & \mapsto (f, !)
\end{align*}
and hence $(\b \eta S_c)_1 \colon (S^21)^\op \times S1 \to \cat{Set}$ and
$(S_c\b \eta )_1 \colon (S^21)^\op \times S1 \to \cat{Set}$ are given by:
\begin{align*}
(\b \eta S_c)_1(\phi; n) &= (\eta_{S1})_\ast(\phi; n) = S^21(\phi, (n \xrightarrow{\id} n))\\
(S_c\b \eta)_1(\phi; n) & = \hat S(\eta_1)_\ast(\phi; n) \cong (S
\eta_1)_\ast(\phi; n) = S^21(\phi, (n \xrightarrow{!} 1))
\end{align*}
Thus the composite along the upper side of this diagram is given by
\[
(\b \eta S_c)_1(\phi; n) = S^21(\phi, (n \xrightarrow{!} 1)) \cong
\begin{cases}
S1(n_\phi, n) & \text{if $m_\phi = 1$;}\\
\emptyset & \text{otherwise,}
\end{cases}\eqno{(1)}
\]
where the isomorphism is natural in $\phi$ and $n$; and with respect to this
isomorphism, the projection down to $S\I$ is given simply by the inclusion
\[(\b \eta S_c)_1(\phi; n) \hookrightarrow S1(n_\phi, n)\text.\]
Now, the lower side is given by
\[
(\b \delta_1 \otimes (S_c \b\eta)_1)(\phi; n) = \int^{\psi \in S^21}
S^21(\psi, (n \xrightarrow{\id} n)) \times \delta_1(\phi; \psi)\text,
\]
which is isomorphic to $\b \delta_1(\phi; (n \xrightarrow{\id} n))$, naturally
in $\phi$ and $n$. Now, any element $f$ of $\b \delta_1(\phi; (n
\xrightarrow{\id} n))$, given by
\[\cd{
  n_\phi \ar[d]_{\phi} \ar[r]^{f_n} & n \ar[d]^{\id} \\
  m_\phi & n
}\] say, must satisfy $m_\phi + n = n_\phi + 1$; but since $n =
n_\phi$, this can only happen if $m_\phi = 1$; and in this case, the diagram
\[\cd{
  n_\phi \ar[d]_{\phi} \ar[r]^{f_n} & n \ar[d]^{!} \\
  m_\phi \ar[r]_{!} & 1
}\] is necessarily a pushout. Hence
\[(\b \delta_1 \otimes (S_c \b\eta)_1)(\phi; n) \cong \begin{cases}
S1(n_\phi, n) & \text{if $m_\phi = 1$;}\\
\emptyset & \text{otherwise,}
\end{cases}\eqno{(2)}
\]
naturally in $\phi$ and $n$; and once again, the projection down to $S\I$ is
given simply by inclusion. So, composing the isomorphisms (1) and (2), we get a
special invertible cell $\overline{\b \eta}_1$ which is compatible with the
projections down to $S\I$, as required. \qedhere
\end{proof}

\begin{Prop}
There is an invertible special cell
\[\cd[@+1em]{
 S_cS_m1
  \rtwocell[0.35]{dr}{\overline {\b \epsilon}_1}
  \ar[r]|-{\object@{|}}^{\b \delta_1}
  \ar[d]|-{\object@{|}}_{(\b \epsilon S_m)_1} &
 S_mS_c1
  \ar[dl]|-{\object@{|}}^{(S_m\b \epsilon)_1} \\
 S_m1 & {}
}\] mediating the centre of this diagram in $\coll(S)$ (where we
omit the projections to $S\I$).
\end{Prop}
\begin{proof}
Dual to the above. \qedhere
\end{proof}

\subsection{(PDD3)}
For (PDD3) we must produce the component of the invertible special
modifications $\overline {\b \mu}$ and $\overline {\b \Delta}$ at $1$:
\begin{Prop}
There is an invertible special cell
\[\cd[@+1em]{
 S_cS_m1
  \dtwocell{drr}{\overline {\b \Delta}_1}
  \ar[d]|-{\object@{|}}_{(\b \Delta S_m)_1}
  \ar[rr]|-{\object@{|}}^{\b \delta_1} & &
 S_mS_c1
  \ar[d]|-{\object@{|}}^{(S_m \b \Delta)_1} \\
 S_cS_cS_m1
  \ar[r]|-{\object@{|}}_{(S_c\b \delta)_1} &
 S_cS_mS_c1
  \ar[r]|-{\object@{|}}_{(\b \delta S_c)_1} &
 S_mS_cS_c1
}\] mediating the centre of this diagram in $\coll(S)$ (where we
omit the projections to $S\I$).
\end{Prop}

\begin{proof}
Let us describe explicitly the horizontal arrows involved in the above diagram.
The functors $\mu_{S1} \colon S^31 \to S^21$ and $S\mu_1 \colon S^31 \to S^21$
in $\Cat$ are given by
\begin{align*}
\mu_{S1} \colon (n_\phi \xrightarrow{\phi_1} m_\phi \xrightarrow{\phi_2}
r_\phi) & \mapsto (n_\phi \xrightarrow{\phi_1} m_\phi) \\
(f_n, f_m, f_r) & \mapsto (f_n, f_m) \\
\text{and } S\mu_1 \colon (n_\phi \xrightarrow{\phi_1} m_\phi \xrightarrow{\phi_2} r_\phi) & \mapsto (n_\phi \xrightarrow{\phi_2 \phi_1} r_\phi) \\
(f_n, f_m, f_r) & \mapsto (f_n, f_r)
\end{align*}
and hence $(\b \Delta S_m)_1 \colon (S^31)^\op \times S^21 \to \cat{Set}$ and
$(S_m\b \Delta)_1 \colon (S^31)^\op \times S^21 \to \cat{Set}$ are given by:
\begin{align*}
(\b \Delta S_m)_1(\phi; \psi) &= (\mu_{S1})^\ast(\phi; \psi) = S^21((n_\phi \xrightarrow{\phi_1} m_\phi), \psi)\\
(S_m\b \Delta)_1(\phi; \psi) &= \hat S(\mu_{1})^\ast(\phi; \psi)
\cong (S\mu_{1})^\ast(\phi; \psi) = S^21((n_\phi
\xrightarrow{\phi_2 \phi_1} r_\phi), \psi)\text.
\end{align*}
We now wish to describe $(\b \delta S_c)_1$ and $(S_c \b \delta)_1$. It's a
straightforward calculation to see that $(\b \delta S_c)_1 \colon (S^31)^\op
\times S^31 \to \cat{Set}$ is given as follows:
\begin{itemize}
\item \textbf{On objects}: elements $f \in (\b \delta S_c)_1(\phi; \psi)$
    are pairs of bijections $f_n$ and $f_m$ fitting in the diagram
\[
\cd[@!@-1em]{
  n_\phi \ar[r]^{f_n} \ar[d]_{\phi_1} & n_\psi \ar[d]^{\psi_1} \\
  m_\phi \ar[r]_{f_m} \ar[d]_{\phi_2} & m_\psi \ar[d]^{\psi_2} \\
  r_\phi & r_\psi
}
\]
such that the span $r_\phi \xleftarrow{\phi_2} m_\phi \xrightarrow{\psi_2
\circ f_m} r_\psi$ is acyclic and connected.
\item \textbf{On maps}: Let $g \colon \psi \to \rho$ in $S^31$ and let $f
    \in (\b \delta S_c)_1(\phi; \psi)$. Then we give an element $g \cdot f
    \in (\b \delta S_c)_1(\phi; \rho)$ by
\[
\cd[@!]{
  n_\phi \ar[r]^{g_n \circ f_n} \ar[d]_{\phi_1} & n_\rho \ar[d]^{\rho_1} \\
  m_\phi \ar[r]_{g_m \circ f_m} \ar[d]_{\phi_2} & m_\rho \ar[d]^{\rho_2} \\
  r_\phi & r_\rho\text;
}
\]
and we give the right action of $S^31$ similarly.
\end{itemize}
Likewise, it's easy to calculate that $(S_c \b \delta)_1 \colon (S^31)^\op
\times S^31 \to \cat{Set}$ is given by:
\begin{itemize}
\item \textbf{On objects}: elements $f \in (S_c \b \delta)_1(\phi; \psi)$
    are pairs of bijections $f_n \colon n_\phi \to n_\psi$ and $f_r \colon
    r_\phi \to r_\psi$ fitting in the diagram
\[
\cd[@!@-1em]{
  n_\phi \ar[r]^{f_n} \ar[d]_{\phi_1} & n_\psi \ar[d]^{\psi_1} \\
  m_\phi \ar[d]_{\phi_2} & m_\psi \ar[d]^{\psi_2} \\
  r_\phi \ar[r]_{f_r} & r_\psi
}
\]
  such that for each $i = 1, \dots, r_\psi$, the induced spans
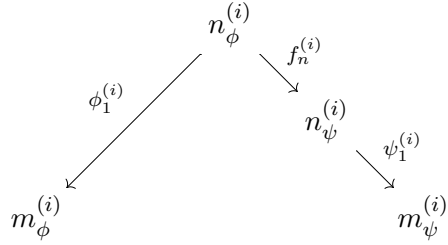
\[\cd[@!@-=1em]{& & n_\phi^{(i)} \ar[dr]^{f_n^{(i)}} \ar[ddll]_{\phi_1^{(i)}} \\
      & & & n_\psi^{(i)} \ar[dr]^{\psi_1^{(i)}} \\
      m_\phi^{(i)} & & & & m_\psi^{(i)}
}
\]
are acyclic and connected.
\end{itemize}
[Let us clarify what the induced spans referred to above actually are. We have
the commuting diagram
\[
\cd[@!@-1em]{
  n_\phi \ar[r]^{f_n} \ar[d]_{\phi_1} & n_\psi \ar[r]^{\psi_1} & m_\psi \ar[d]^{\psi_2} \\
  m_\phi \ar[r]_{\phi_2} & r_\phi \ar[r]_{f_r} & r_\psi
}\eqno{(*)}
\]
and the induced spans are the result of pulling this diagram back along
elements $i \colon 1 \to r_\psi$. By the results of the first section of this
chapter, these spans are all acyclic and connected if and only if $(*)$ is a
pushout and $r_\psi + n_\phi = m_\phi + m_\psi$.]
\begin{itemize}
\item \textbf{On maps}: Let $g \colon \psi \to \rho$ in $S^31$ and let $f
    \in (S_c \b \delta)_1(\phi; \psi)$. Then we give an element $g \cdot f
    \in (S_c \b \delta)_1(\phi; \rho)$ by
\[
\cd[@!]{
  n_\phi \ar[r]^{g_n \circ f_n} \ar[d]_{\phi_1} & n_\rho \ar[d]^{\rho_1} \\
  m_\phi \ar[d]_{\phi_2} & m_\rho \ar[d]^{\rho_2} \\
  r_\phi \ar[r]_{g_r \circ f_r} & r_\rho\text;
}
\]
and we give the right action similarly.
\end{itemize}
Now, returning to the diagram in question, the upper side is given by
\[
((S_m \b \Delta)_1 \otimes \b \delta_1)(\phi; \rho) = \int^{\psi \in S^21}
\b \delta_1(\psi; \rho) \times S^21\big((n_\phi \xrightarrow{\phi_2 \phi_1}
r_\phi), \psi\big)\text,
\]
which is isomorphic to $\b \delta_1((n_\phi \xrightarrow{\phi_2 \phi_1}
r_\phi); \rho)$, naturally in $\phi$ and $\rho$. With respect to this
isomorphism, the projection onto $S\I$ has component morphisms $\b
\delta_1\big((n_\phi \xrightarrow{\phi_2 \phi_1} r_\phi); \rho\big) \to
S1(n_\phi; n_\rho)$ which send
\[
\cd[@!@-1em]{
  n_\phi \ar[r]^{f_n} \ar[d]_{\phi_2 \phi_1} & n_\rho \ar[d]^\rho \\
  r_\phi  & m_\rho
}
\]
to $f_n$. The lower side of this diagram, which we denote by $K$, is given by
\begin{align*}
 K(\phi; \rho) &= \big((\b \delta S_c)_1 \otimes (S_c \b \delta)_1 \otimes (\b \Delta
S_m)_1\big)(\phi; \rho)\\ & = \int^{\psi, \xi \in S^31} S^21((n_{\xi} \xrightarrow{\xi_1} m_{\xi}), \rho) \times
(S_c \b \delta)_1(\psi; \xi) \times (\b \delta S_c)_1(\phi; \psi)\text.
\end{align*}
We may represent a typical element $x \in K(\phi; \rho)$ as $x = f \otimes g
\otimes h$, where $f \in (\b \delta S_c)_1(\phi; \psi)$, $g \in (S_c \b
\delta)_1(\psi; \xi)$, and $h \in S^21((n_{\xi} \xrightarrow{\xi_1} m_{\xi}),
\rho)$:
\[
\cd[@!@-1em]{
  n_\phi \ar[r]^{f_n} \ar[d]_{\phi_1} & n_\psi \ar[r]^{g_n} \ar[d]_{\psi_1} & n_\xi \ar[r]^{h_n} \ar[d]_{\xi_1} & n_\rho \ar[d]^{\rho} \\
  m_\phi \ar[r]_{f_m} \ar[d]_{\phi_2} & m_\psi \ar[d]_{\psi_2} & m_\xi \ar[r]_{h_m} \ar[d]_{\xi_2} & m_\rho\\
  r_\phi  & r_\psi \ar[r]_{g_r} & r_\xi\text.
}
\]
Then the projection onto $S\I$ has components
\begin{align*}
K(\phi; \rho) & \to S1(n_\phi, n_\rho)\\
f \otimes g \otimes h  & \mapsto h_n \circ g_n \circ f_n\text.
\end{align*}
So, we need to set up an isomorphism between $K(\phi; \rho)$ and $\b
\delta_1((n_\phi \xrightarrow{\phi_2 \phi_1} r_\phi); \rho)$ which is natural
in $\phi$ and $\rho$ and compatible with the projections onto $S\I$. In one
direction, we send the element $x \in K(\phi; \rho)$:
\[
\cd[@!@-1em]{
  n_\phi \ar[r]^{f_n} \ar[d]_{\phi_1} & n_\psi \ar[r]^{g_n} \ar[d]_{\psi_1} & n_\xi \ar[r]^{h_n} \ar[d]_{\xi_1} & n_\rho \ar[d]^{\rho} \\
  m_\phi \ar[r]_{f_m} \ar[d]_{\phi_2} & m_\psi \ar[d]_{\psi_2} & m_\xi \ar[r]_{h_m} \ar[d]_{\xi_2} & m_\rho\\
  r_\phi  & r_\psi \ar[r]_{g_r} & r_\xi
}
\]
to the element $\hat x$ of $\b \delta_1\big((n_\phi \xrightarrow{\phi_2 \phi_1}
r_\phi); \rho\big)$ given by
\[
\cd[@!]{
  n_\phi \ar[r]^{h_n g_n f_n} \ar[d]_{\phi_2 \phi_1} & n_\rho \ar[d]^{\rho} \\
  r_\phi & m_\rho\text.
}
\]
Note that this element is independent of the representation of $x$ that we
chose, that this assignation is natural in $\phi$ and $\rho$, and is compatible
with the projection down to $S\I$; but for it to be well-defined, we need still
to check that the span $r_\phi \xleftarrow{\phi_2 \phi_1} n_\phi
\xrightarrow{\rho h_n g_n f_n} m_\rho$ is acyclic and connected. For this, we
observe first that in the following diagram
\[
\cd[@!@-1em]{
  n_\phi \ar[r]^{f_n} \ar[d]_{\phi_1} & n_\psi \ar[r]^{g_n} \ar[d]_{\psi_1} & n_\xi \ar[r]^{\xi_1} & m_\xi \ar[d]^{\xi_2} \ar[r]^{h_n} & m_\rho \ar[dd]\\
  m_\phi \ar[d]_{\phi_2} \ar[r]_{f_m} & m_\psi \ar[r]_{\psi_2} & r_\psi \ar[r]_{g_r} \ar[d] & r_\xi \ar[d]\\
  r_\phi \ar[rr] & & 1 \ar[r] & 1 \ar[r] & 1
}
\]
each of the smaller squares is a pushout; and hence the outer square is also a
pushout. But the top edge is $h_n \xi_1 g_n f_n = \rho h_n g_n f_n$, so that
the square
\[
\cd[@!@+1em]{
  n_\phi \ar[r]^{\rho h_n g_n f_n} \ar[d]_{\phi_2 \phi 1} & n_\rho \ar[d] \\
  r_\psi \ar[r] & 1
}
\]
is a pushout as required. Furthermore, the following equalities hold:
\begin{align*}
  r_\phi + r_\psi &= m_\phi + 1\text, &
  m_\psi + m_\xi &= n_\psi + r_\xi\text,\\
  m_\psi &= m_\phi\text, &
  m_\rho &= m_\xi\text,\\
  r_\psi &= r_\xi\text, &
  \text{and \ \ \ } n_\psi &= n_\phi
\end{align*}
whence we have $m_\rho + r_\phi = n_\phi + 1$. So the span $r_\phi
\xleftarrow{\phi_2 \phi_1} n_\phi \xrightarrow{\rho h_n g_n f_n} m_\rho$ is
acyclic and connected as required.

Conversely, suppose we are given an element $k$ of $\b \delta_1((n_\phi
\xrightarrow{\phi_2 \phi_1} r_\phi); \rho)$:
\[
\cd[@!]{
  n_\phi \ar[r]^{k_n} \ar[d]_{\phi_2 \phi_1} & n_\rho \ar[d]^{\rho} \\
  r_\phi & m_\rho\text;
}
\]
then we take the following pushout:
\[
\cd[@!]{
  n_\phi \ar[r]^{\rho k_n} \ar[d]_{\phi_1} & m_\rho \ar[d]^{i_2} \\
  m_\phi \ar[r]_{i_1} & r\text.
}
\]
Now, the map $i_1$ in this pushout square need not be order-preserving; but it
has a (non-unique) factorisation as $m_\phi \xrightarrow{\alpha_1} r_1
\xrightarrow{\sigma_1} r$, where $\alpha_1$ is order-preserving and $\sigma_1$
a bijection. Similarly, we can factorise $i_2$ as $m_\rho
\xrightarrow{\alpha_2} r_2 \xrightarrow{\sigma_2} r$ with $\alpha_2$ is
order-preserving and $\sigma_2$ a bijection. [Note that it follows that each of
the diagrams
\[
\cd[@!]{
  n_\phi \ar[r]^{\rho k_n} \ar[d]_{\phi_1} & m_\rho \ar[d]^{\sigma_1^{-1} i_2} \\
  m_\phi \ar[r]_{\alpha_1} & r_1
}\quad \text{and} \quad \cd[@!]{
  n_\phi \ar[r]^{\rho k_n} \ar[d]_{\phi_1} & m_\rho \ar[d]^{\alpha_2} \\
  m_\phi \ar[r]_{\sigma_2^{-1} i_1} & r_2
}
\]
is also a pushout.] Now we send $k$ to the element $\hat k$ of $K(\phi; \rho)$
represented by the following:
\[
\cd[@!@-1em]{
  n_\phi \ar[r]^{\id} \ar[d]_{\phi_1} & n_\phi \ar[r]^{k_n} \ar[d]_{\phi_1} & n_\rho \ar[r]^{\id} \ar[d]_{\rho} & n_\rho \ar[d]^{\rho} \\
  m_\phi \ar[r]_{\id} \ar[d]_{\phi_2} & m_\phi \ar[d]_{\alpha_1} & m_\rho \ar[r]_{\id} \ar[d]_{\alpha_2} & m_\rho\text.\\
  r_\phi  & r_1 \ar[r]_{\sigma_2^{-1} \sigma_1} & r_2
}
\]
This is visibly compatible with the projection down onto $S\I$, but we need to
check that it is in fact a valid element of $K(\phi; \rho)$. Clearly all
squares commute in the diagram above, so we need only check the acyclic and
connected conditions. We start with connectedness; for the middle map, the
diagram
\[
\cd[@!]{
  n_\phi \ar[r]^{k_n} \ar[d]_{\phi_1} & n_\rho \ar[r]^{\rho} & m_\rho \ar[d]^{\alpha_2}\\
  m_\phi \ar[r]_{\alpha_1} & r_1 \ar[r]_{\sigma_2^{-1} \sigma_1} & r_2
} = \cd[@!]{
  n_\phi \ar[r]^{\rho k_n} \ar[d]_{\phi_1} & m_\rho \ar[d]^{\alpha_2} \\
  m_\phi \ar[r]_{\sigma_2^{-1} i_1} & r_2
}
\]
is indeed a pushout, so the induced spans for the middle map are connected. For
the left-hand map, consider the diagram
\[
\cd{
  n_\phi \ar[r]^{\rho k_n} \ar[d]_{\phi_1} & m_\rho \ar[d]^{\sigma_1^{-1} i_2} \\
  m_\phi \ar[r]_{\alpha_1} \ar[d]_{\phi_2} & r_1 \ar[d]\\
  r_\phi \ar[r] & 1\text;
}
\]
the outer square and the upper square are both pushouts, and hence so is the
lower square; so the left-hand span is connected.

And now acyclicity. For the middle map, we need that, given any monomorphism
$\iota \colon n_\phi' \hookrightarrow n_\phi$, the diagram
\[
\cd{
  n_\phi' \ar[r]^{\rho k_n \iota} \ar[d]_{\phi_1  \iota} & m_\rho \ar[d]^{\alpha_2} \\
  m_\phi \ar[r]_{\sigma_2^{-1} i_1} & r_2
}
\]
is no longer a pushout. But suppose it were; then in the diagram
\[
\cd{
  n_\phi' \ar[r]^{\rho k_n \iota} \ar[d]_{\phi_1 \iota} & m_\rho \ar[d]^{\sigma_1^{-1} i_2} \\
  m_\phi \ar[r]_{\alpha_1} \ar[d]_{\phi_2} & r_1 \ar[d]\\
  r_\phi \ar[r] & 1
}
\]
the upper and lower squares would be pushouts, hence making the outer edge a
pushout; but this contradicts the acyclicity of the span $r_\phi \leftarrow
n_\phi \rightarrow m_\rho$. So the induced spans for the middle map are
acyclic. Thus we now know that the following equations hold:
\begin{align*}
m_\phi + m_\rho &= n_\phi + r_2 \\
r_\phi + m_\rho &= n_\phi + 1\\
r_1 &= r_2\text,
\end{align*}
and so can deduce that $r_1 + r_\phi = m_\phi + 1$, as required for the
left-hand span to be acyclic.

It remains to check that these two assignations are mutually inverse. It is
evident, given $k \in d_1((n_\phi \xrightarrow{\phi_2 \phi_1} r_\phi); \rho)$,
that $\hat{\hat k} = k$. For the other direction, we send
\[
x = \cd[@!@-1em]{
  n_\phi \ar[r]^{f_n} \ar[d]_{\phi_1} & n_\psi \ar[r]^{g_n} \ar[d]_{\psi_1} & n_\xi \ar[r]^{h_n} \ar[d]_{\xi_1} & n_\rho \ar[d]^{\rho} \\
  m_\phi \ar[r]_{f_m} \ar[d]_{\phi_2} & m_\psi \ar[d]_{\psi_2} & m_\xi \ar[r]_{h_m} \ar[d]_{\xi_2} & m_\rho\\
  r_\phi  & r_\psi \ar[r]_{g_r} & r_\xi\text.
}\quad \text{to} \quad \hat{\hat x} = \cd[@!@-1em]{
  n_\phi \ar[r]^{\id} \ar[d]_{\phi_1} & n_\phi \ar[r]^{k_n} \ar[d]_{\phi_1} & n_\rho \ar[r]^{\id} \ar[d]_{\rho} & n_\rho \ar[d]^{\rho} \\
  m_\phi \ar[r]_{\id} \ar[d]_{\phi_2} & m_\phi \ar[d]_{\alpha_1} & m_\rho \ar[r]_{\id} \ar[d]_{\alpha_2} & m_\rho\\
  r_\phi  & r_1 \ar[r]_{\sigma_2^{-1} \sigma_1} & r_2\text.
}
\]
We claim that these two diagrams represent the same element of $K(\phi; \rho)$.
Indeed, note that in the diagram
\[\cd[@!@-1em]{
  n_\phi \ar[r]^{f_n} \ar[d]_{\phi_1} & n_\psi \ar[r]^{g_n} \ar[d]_{\psi_1} &
  n_\xi \ar[r]^{\xi_1} & m_\xi \ar[r]^{h_m} \ar[d]_{\xi_2} & m_\rho \ar[d]^{g_r^{-1}\xi_2 h_m^{-1}}\\
  m_\phi \ar[r]_{f_m} & m_\psi \ar[r]_{\psi_2} & r_\psi \ar[r]_{g_r} & r_\xi \ar[r]_{g_r^{-1}} & r_\psi}\]
each of the smaller squares is a pushout, and hence the outer edge is. But the
upper edge is $h_m \xi_1 g_n f_n = \rho h_n g_n f_n = \rho k_n$, so that the
diagram
\[\cd[@!@-1em]{
  n_\phi \ar[r]^{\rho k_n} \ar[d]_{\phi_1} & m_\rho \ar[d]^{g_r^{-1}\xi_2 h_m^{-1}}\\
  m_\phi \ar[r]_{\psi_2 f_m} & r_\psi}\] is a pushout. Since $r_1$ is
also a pushout for this diagram, it follows that there is an isomorphism
$\beta_1 \colon r_1 \to r_\psi$ such that $\beta_1 \alpha_1 = \psi_2 f_m$;
hence the following diagram commutes:
\[
\cd[@!@-1em]{
n_\phi \ar[d]_{\phi_1} \ar[r]^{f_n} & n_\psi \ar[d]^{\psi_1} \\
m_\phi \ar[d]_{\alpha_1} \ar[r]^{f_m} & m_\psi \ar[d]^{\psi_2} \\
r_1 \ar[r]_{\beta_1} & r_\psi }
\]
Similarly, we see that
\[\cd[@!@-1em]{
  n_\phi \ar[r]^{\rho k_n} \ar[d]_{\phi_1} & m_\rho \ar[d]^{\xi_2 h_m^{-1}}\\
  m_\phi \ar[r]_{g_r \psi_2 f_m} & r_\xi}\] is a pushout, and so
there is an isomorphism $\beta_2 \colon r_\xi \to r_2$ such that $\beta_2 \xi_2
h_m^{-1} = \alpha_2$, i.e., $\beta_2 \xi_2 = \alpha_2 h_m$. Hence the following
diagram commutes:
\[
\cd[@!@-1em]{
n_\xi \ar[d]_{\xi_1} \ar[r]^{h_n} & n_\rho \ar[d]^{\rho} \\
m_\xi \ar[d]_{\xi_2} \ar[r]^{h_m} & m_\rho \ar[d]^{\alpha_2} \\
r_\xi \ar[r]_{\beta_2} & r_2\text. }
\]
Furthermore, we have $r_1 \xrightarrow{\beta_1} r_\psi \xrightarrow{g_r} r_\xi
\xrightarrow{\beta_2} r_2 = r_1 \xrightarrow{\sigma_1} r
\xrightarrow{\sigma_2^{-1}} r_2$, since each of these objects is a pushout of
the same span, and the isomorphisms between them are isomorphisms of pushouts.
Thus, using an evident notation for the internal actions, we have
\begin{align*}
x = \cd[@!@-1em]{
  n_\phi \ar[r]^{f_n} \ar[d]_{\phi_1} & n_\psi \ar[r]^{g_n} \ar[d]_{\psi_1} & n_\xi \ar[r]^{h_n} \ar[d]_{\xi_1} & n_\rho \ar[d]^{\rho} \\
  m_\phi \ar[r]_{f_m} \ar[d]_{\phi_2} & m_\psi \ar[d]_{\psi_2} & m_\xi \ar[r]_{h_m} \ar[d]_{\xi_2} & m_\rho\\
  r_\phi  & r_\psi \ar[r]_{g_r} & r_\xi\text.
}& \equiv \cd[@!@-1em]{
  n_\phi \ar[r]^{\id} \ar[d]_{\phi_1} & n_\phi \ar[d]_{\phi_1} \ar[r]^{f_n} & n_\psi \ar[d]^{\psi_1} \ar[r]^{g_n} &
  n_\xi \ar[d]_{\xi_1} \ar[r]^{h_n} & n_\rho \ar[r]^{\id} \ar[d]_{\rho} & n_\rho \ar[d]^{\rho} \\
  m_\phi \ar[r]_{\id} \ar[d]_{\phi_2} &
  m_\phi \ar[d]_{\alpha_1} \ar[r]^{f_m} & m_\psi \ar[d]^{\psi_2} &
  m_\xi \ar[d]_{\xi_2} \ar[r]^{h_m} & m_\rho \ar[r]_{\id} \ar[d]_{\alpha_2} & m_\rho\\
  r_\phi  & r_1 \ar[r]_{\beta_1} & r_\psi \ar[r]_{g_r} & r_\xi \ar[r]_{\beta_2} & r_2\text.
}\\
& \equiv \cd[@!@-1em]{
  n_\phi \ar[r]^{\id} \ar[d]_{\phi_1} & n_\phi \ar[r]^{k_n} \ar[d]_{\phi_1} & n_\rho \ar[r]^{\id} \ar[d]_{\rho} & n_\rho \ar[d]^{\rho} \\
  m_\phi \ar[r]_{\id} \ar[d]_{\phi_2} & m_\phi \ar[d]_{\alpha_1} & m_\rho \ar[r]_{\id} \ar[d]_{\alpha_2} & m_\rho\\
  r_\phi  & r_1 \ar[r]_{\sigma_2^{-1} \sigma_1} & r_2\text.
} = \hat{\hat x}\text.
\end{align*}

So the assignations $x \mapsto \hat x$ and $k \mapsto \hat k$ are mutually
inverse as required. It now follows that the assignation $\b \delta_1((n_\phi
\xrightarrow{\phi_2 \phi_1} r_\phi); \rho) \to K(\phi; \rho)$ is natural in
$\phi$ and $\rho$, since its inverse is. \qedhere
\end{proof}

\begin{Prop}
There is an invertible special cell
\[
\cd[@+1em]{
 S_cS_mS_m1
  \dtwocell{drr}{\overline {\b \mu}_1}
  \ar[r]|-{\object@{|}}^{(\b \delta S_m)_1}
  \ar[d]|-{\object@{|}}_{(S_c\b \mu)_1} &
 S_mS_cS_m1
  \ar[r]|-{\object@{|}}^{(S_m\b \delta)_1} &
 S_cS_cS_m1
  \ar[d]|-{\object@{|}}^{(\b \mu S_c)_1} \\
 S_cS_m1
  \ar[rr]|-{\object@{|}}_{\b \delta_1} & &
 S_mS_c1
}
\]
mediating the centre of this diagram in $\coll(S)$ (where we omit the
projections to $S\I$).
\end{Prop}
\begin{proof}
Dual to the above. \qedhere
\end{proof}

\subsection{(PDA1)--(PDA10)}
It remains only to show that the data produced above satisfies the ten
coherence axioms  (PDA1)--(PDA10). At first this may appear somewhat
forbidding, but our job is made rather simple by the following argument.

\begin{Defn}
We say that a cell
\[\cd{
X_s \ar[d]_{f_s} \ar[r]|-{\object@{|}}^{\X} \dtwocell{dr}{\b f}& X_t \ar[d]^{f_t}  \\
Y_s \ar[r]|-{\object@{|}}_{\Y} & Y_t }\] of $\dCat$ is
\defn{locally monomorphic} if it is a monomorphism when viewed as
a map of $[X_t^\op \times X_s, \cat{Set}]$:
\[\cd{
 X_t^\op \times X_s
  \ar[rr]^{f_t^\op \times f_s}
  \ar[dr]_{X} &
 {}
  \rtwocell{d}{f} &
 Y_t^\op \times Y_s
  \ar[dl]^Y \\ &
 \cat{Set}\text.
}\]
\end{Defn}

Now, local monomorphisms admit a limited form of `left cancellation'. Indeed,
suppose we are given objects $\X = X \colon X_s \tor X_t$ and $\X' = X' \colon
X_s \tor X_t$ of $\dCat_1$, and \emph{special} maps $\b g_1$ and $\b g_2 \colon
\X' \to \X$; then given a local monomorphism $\b f \colon \X \to \Y$, we have
that
\[\b f \circ \b g_1 = \b f \circ \b g_2 \quad \text{implies} \quad
\b g_1 = \b g_2\text,\] since to give a special map $g_i \colon
\X' \to \X$ is equivalently to give a natural transformation $g_i \colon X'
\Rightarrow X$; therefore the result follows from the fact that $f \colon X
\Rightarrow (Y \circ f_t^\op \times f_s)$ is a monomorphism in $[X_t^\op \times
X_s, \cat{Set}]$.

Observe also that, given a special isomorphism $\b g \colon \X' \to \X$ and a
local monomorphism $\b f \colon \b X \to \b Y$, the map $\b f \circ \b g$ is
again a local monomorphism.

\begin{Prop}
Consider each of the pasting diagrams in the axioms (PDA1)--(PDA10) as a
diagram in $\dCat / S\I_1$. Then the projection map from each `source' and
`target' face down onto $S\I_1$ is a local monomorphism.
\end{Prop}
\begin{proof}
Observe that every special cell in the pasting diagrams for (PDA1)--(PDA10) is
invertible, and therefore, for each pasting diagram it suffices to show for
\emph{any one} path through it that the projection onto $S\I_1$ is a local
monomorphism; it then follows, by the discussion preceding this proposition,
that the same is true for all other paths. We now work our way through the ten
axioms:
\begin{itemize}
\item{(PDA1):} Let us write $K$ for the composite $S_c1 \xrightarrow{\b
    \epsilon_1} \id \xrightarrow{\b \eta_1} S_m1$; then we have
\[K(m; n) = \begin{cases} \{\,\ast\,\} & \text{if $m = n = 1$;} \\ \emptyset &
\text{otherwise.}\end{cases}\] and the projection down onto
$S\I_1$ simply sends the unique element of $K(1; 1)$ to the unique element
of $S1(1; 1)$, and thus is a local monomorphism as required.
\item{(PDA2)--(PDA5):} For each of these we look at the path $\b \delta_1
    \colon S_cS_m1 \to S_mS_c1$, and from the definitions, the projection
    onto $S\I_1$ is visibly a local monomorphism.
\item{(PDA6):} Let us write $K$ for the composite \[S_cS_mS_mS_m1
    \xrightarrow{(S_cS_m \b \mu)_1} S_cS_mS_m1 \xrightarrow{(S_c \b \mu)_1}
    S_cS_m1 \xrightarrow{\b \delta_1} S_mS_c1\text.\] Then we have an
    isomorphism
\[K(\phi; \psi) \cong \b \delta_1\big(\phi; (n_\psi \xrightarrow{\psi_3\psi_2\psi_1} s_\psi)\big)\]
natural in $\phi$ and $\psi$, where we are writing a typical element of
$S_cS_mS_mS_m1$ as $\psi = n_\psi \xrightarrow{\psi_1} m_\psi
\xrightarrow{\psi_2} r_\psi \xrightarrow{\psi_3} s_\psi$ in the evident
way. With respect to this isomorphism, the projection down onto $S\I_1$ is
given simply by the value of $\tilde{\b \delta}_1$ there, which is a
monomorphism as required.
\item{(PDA7):} Dual to (PDA6).
\item{(PDA8):} Let us write $K$ for the composite
\begin{equation*}
    S_cS_mS_m1 \xrightarrow{(S_c \b
    \mu)_1} S_cS_m1 \xrightarrow{\b \delta_1} S_mS_c1 \xrightarrow{(S_m \b \epsilon)_1} S_m1\text;
\end{equation*}
then we have
\[K(m; \phi) \cong \b \delta_1\big((m \xrightarrow{\id} m); (n_\phi \xrightarrow{\phi_2 \phi_1} r_\phi)\big)\]
and again the projection down onto $S\I_1$ is simply given by the value of
$\tilde {\b \delta}_1$ there; and so a local monomorphism.
\item{(PDA9):} Dual to (PDA8).
\item{(PDA10):} Let us write $K$ for the composite
\begin{equation*}
    S_cS_mS_m1 \xrightarrow{(S_c \b \mu)_1} S_cS_m1 \xrightarrow{\b \delta_1} S_cS_m1 \xrightarrow{(S_m
    \b \Delta)_1} S_mS_cS_c1\text;
\end{equation*}
then we have
\[K(\psi; \phi) \cong \b \delta_1\big((n_\psi \xrightarrow{\psi_2 \psi_1} r_\psi); (n_\phi \xrightarrow{\phi_2 \phi_1} r_\phi)\big)\text.\]
Once more, the projection down onto $S\I_1$ is just the value of $\tilde
{\b \delta}_1$ there, and so a local monomorphism. \qedhere
\end{itemize}
\end{proof}

\begin{Cor}
The pasting equalities (PDA1)--(PDA10), when viewed as diagrams in $\dCat /
S\I_1$,  hold for the data (PDD1)--(PDD5) given above.
\end{Cor}
\begin{proof}
Consider (PDA1) for example. The two pasting diagrams under consideration pick
out two arrows $\b f$ and $\b g$ of $\dCat_1 / S\I_1$:
\[\cd{
  (\b \epsilon S_m)_1 \otimes (S_c \b \eta)_1
   \ar[rr]^-{\b f}
   \ar[dr]_{\b \pi_1} & &
  (S_m \b \epsilon)_1 \otimes (\b \eta S_c)_1
   \ar[dl]^{\b \pi_2} \\ &
  S\I_1
}\] and
\[\cd{
  (\b \epsilon S_m)_1 \otimes (S_c \b \eta)_1
   \ar[rr]^-{\b g}
   \ar[dr]_{\b \pi_1} & &
  (S_m \b \epsilon)_1 \otimes (\b \eta S_c)_1
   \ar[dl]^{\b \pi_2} \\ &
  S\I_1\text,
}
\]
where both the above diagrams commute. But by the previous proposition, the
projections $\b \pi_1$ and $\b \pi_2$ are local monomorphisms, and since $\b f$
and $\b g$ are special maps, we have
\[\b \pi_2 \circ \b f = \b \pi_1 = \b \pi_2 \circ \b g \quad \text{implying} \quad \b f = \b g\text.\]
We argue similarly for the other nine diagrams. \qedhere
\end{proof}

\noindent This completes the definition of our pseudo-distributive law in
$\mathcal B(\dCat / S\I_1)$; so now, by the arguments of Section 2, we can
produce from this a pseudo-distributive law in $\mathcal B\big(\coll(S)\big)$,
and thence, via the strict homomorphism $V \colon \mathcal B\big(\coll(S)\big)
\to [\cat{Mod}, \cat{Mod}]_\psi$, our desired pseudo-distributive law $\delta
\colon \hat S_c \hat S_m \Rightarrow \hat S_m \hat S_c$ in $\cat{Mod}$.

We are now finally able to state our abstract description of polycategories:

\begin{Defn}
A \defn{polycompositional polycategory} with object set $X$ is a monad on the
discrete object $X$ in the bicategory $Kl(\delta)$.
\end{Defn}

There is one loose end to tie up: we must complete the argument begun in
Proposition \ref{full}, and show that the polycompositional polycategories we
have just defined are equivalent to polycategories equipped with a binary
composition.

\begin{Prop}\label{characterise}
There is a bijection between polycompositional polycategories with object set
$X$; and polycategories with object set $X$ in the sense of Definition
\ref{defpoly}.
\end{Prop}
\begin{proof}
From the arguments which conclude Section \ref{polycatssec}, together with the
explicit description of $\delta_X$ given at the end of \S \ref{pdd1sec}, we see
that the basic data for a polycompositional polycategory with object set $X$
are: sets of polymaps, equipped with actions by the symmetric groups; identity
maps $x \to x$ for each element $x \in X$; and polycomposites for each pair of
families of polymaps equipped with a suitable matching.

The axioms which a polycompositional category will satisfy are associativity
and unitality laws, which may be extracted from the axioms for the
corresponding monad in $Kl(\delta)$; and compatibility laws between
polycomposition and exchange isomorphisms, which may be deduced from an
examination of the coend composition in $Kl(\delta)$.

It thus follows from Proposition \ref{full} that we may derive the basic data
for a polycompositional polycategory from the data for a standard polycategory,
and vice versa; and it is now a matter of straightforward verification to check
that the axioms for the one entail the axioms for the other. Thus we have
assignations in both directions between standard polycategories to
polycompositional polycategories; and further verification shows these
assignations to be mutually inverse. \qedhere
\end{proof}

And so we conclude with the main result of this paper:

\begin{Thm}
To give a polycategory with object set $X$ is to give a monad on the discrete
object $X$ in the bicategory $Kl(\delta)$.
\end{Thm}

\section*{Appendix: Pseudo notions}
We give here definitions of \emph{pseudomonad}, \emph{pseudocomonad} and of a
\emph{pseudo-distributive law} of the latter over the former.

\begin{Defn}
A \defn{pseudomonad} on a bicategory $\mathcal B$ consists of the following
data:
\begin{enumerate}[(PMD1)]
\item A homomorphism $S \colon \mathcal B \to \mathcal B$;
\item Pseudonatural transformations $\eta \colon \id_{\mathcal B}
    \Rightarrow S$ and $\mu \colon SS \Rightarrow S$;
\item Invertible modifications
\[
\cd[@+1em]{
  S
  \ar[d]_{ S  \eta}
  \ar[dr]^{\id_S} & {} \\
  S  S
  \rtwocell[0.35]{ur}{\lambda}
  \ar[r]_{ \mu} &
  S\text,
} \qquad \cd[@+1em]{
  S
  \ar[d]_{  \eta S}
  \ar[dr]^{\id_S} & {} \\
  S  S
  \rtwocell[0.35]{ur}{\rho}
  \ar[r]_{ \mu} &
  S\text,
} \quad \text{and} \quad \cd[@+1em]{
  S  S  S
  \ar[r]^{ S  \mu}
  \ar[d]_{ \mu  S}
  \rtwocell{dr}{\tau} &
  S  S
  \ar[d]^{ \mu} \\
  S  S
  \ar[r]_{ \mu} &
  S\text.
}\]
\end{enumerate}
All subject to the following two axioms:
\begin{enumerate}[(PM{A}1)]
\item The following pastings agree:
\[\cd{
 S^4
  \ar[rr]^{SS\mu}
  \ar[dr]^{S\mu S}
  \ar[dd]_{\mu SS} & &
 S^3
  \ar[dr]^{S \mu}
  \rtwocell{dl}{S \tau} \\ &
 S^3
  \ar[rr]_{S \mu}
  \ar[dd]_{\mu S}
  \rtwocell{ddrr}{\tau} & &
 S^2
  \ar[dd]^{\mu} \\
 S^3
  \ar[dr]_{\mu S}
  \rtwocell{ur}{\tau S} \\ &
 S^2
  \ar[rr]_{\mu} & &
 S
} \quad = \quad \cd{
 S^4
  \ar[rr]^{SS\mu}
  \rtwocell{ddrr}{\cong}
  \ar[dd]_{\mu SS} & &
 S^3
  \ar[dr]^{S \mu}
  \ar[dd]^{\mu S}
  \\ & & &
 S^2
  \ar[dd]^{\mu} \\
 S^3
  \ar[dr]_{\mu S}
  \ar[rr]_{S \mu} & &
 S^2
  \ar[dr]_{\mu}
  \rtwocell{ur}{\tau}
  \rtwocell{dl}{\tau}
  \\ &
 S^2
  \ar[rr]_{\mu} & &
 S\text;
}\]
\item The following pastings agree:
\[
\cd{ &
 S^3
  \ar[dr]^{S \mu}
  \ar[rr]^{\mu S} & &
 S^2
  \ar[dr]^{\mu} \\
 S^2
  \ar[ur]^{S \eta S}
  \ar[rr]_{\id} &
 {}
  \dtwocell[0.4]{u}{S \rho} &
 S^2
  \ar[rr]_{\mu}
  \dtwocell{ur}{\tau} & &
 S
 } \quad = \quad
\cd{ &
 S^3
  \ar[dr]^{\mu S}\\
 S^2
  \ar[ur]^{S \eta S}
  \ar[rr]_{\id} &
 {}
  \dtwocell[0.4]{u}{\lambda S} &
 S^2
  \ar[rr]_{\mu} & &
 S\text.
 }
\]
\end{enumerate}
\end{Defn}

\noindent Dually, we have the notion of a pseudocomonad on a bicategory:

\begin{Defn}
A \defn{pseudocomonad} on a bicategory $\mathcal B$ consists of the following
data:
\begin{enumerate}[(PCD1)]
\item A homomorphism $T \colon \mathcal B \to \mathcal B$;
\item Pseudonatural transformations $\epsilon \colon T \Rightarrow
    \id_{\mathcal B}$ and $\Delta \colon T \Rightarrow TT$;
\item Invertible modifications
\[
\cd[@+1em]{
  T
  \ar[r]^{\Delta}
  \ar[dr]_{\id_T} &
  T^2\text,
  \rtwocell[0.35]{dl}{\lambda'}
  \ar[d]^{T \epsilon} \\ &
  T
} \qquad \cd[@+1em]{
  T
  \ar[r]^{ \Delta}
  \ar[dr]_{\id_T} &
  T^2
  \rtwocell[0.35]{dl}{\rho'}
  \ar[d]^{ \epsilon T} \\ &
  T
} \quad \text{and} \quad \cd[@+1em]{
  T
  \ar[r]^-{\Delta}
  \ar[d]_{\Delta}
  \rtwocell{dr}{\tau'} &
  T^2
  \ar[d]^{\Delta T} \\
  T^2
  \ar[r]_-{T \Delta} &
  T^3\text.
}\]
\end{enumerate}
Subject to the two axioms:
\begin{enumerate}[(PC{A}1)]
\item The following pastings agree:
\[\cd{
 T
  \ar[rr]^{\Delta}
  \ar[dr]^{\Delta}
  \ar[dd]_{\Delta} & &
 T^2
  \ar[dr]^{\Delta T}
  \rtwocell{dl}{\tau'} \\ &
 T^2
  \ar[rr]_{T \Delta}
  \ar[dd]_{\Delta T}
  \rtwocell{ddrr}{\cong} & &
 T^3
  \ar[dd]^{\Delta TT} \\
 T^2
  \ar[dr]_{T \Delta}
  \rtwocell{ur}{\tau'} \\ &
 T^3
  \ar[rr]_{TT\Delta} & &
 T^4
} \quad = \quad \cd{
 T
  \ar[rr]^{\Delta}
  \rtwocell{ddrr}{\tau'}
  \ar[dd]_{\Delta} & &
 T^2
  \ar[dr]^{\Delta T}
  \ar[dd]^{\Delta T}
  \\ & & &
 T^3
  \ar[dd]^{\Delta TT} \\
 T^2
  \ar[dr]_{T \Delta}
  \ar[rr]_{T \Delta} & &
 T^3
  \ar[dr]_{T \Delta T}
  \rtwocell{ur}{\tau' T}
  \rtwocell{dl}{T\tau'}
  \\ &
 T^3
  \ar[rr]_{TT \Delta} & &
 T^4\text;
}\]
\item The following pastings agree:
\[
\cd{
 T
  \ar[dr]_{\Delta}
  \ar[rr]^{\Delta} & &
 T^2
  \ar[rr]^{\id_{T^2}}
  \ar[dr]_{T\Delta}
  \dtwocell{dl}{\tau'} &
 {}
  \dtwocell{d}{T \rho'} &
 T^2 \\ &
 T^2
  \ar[rr]_{\Delta T} & &
 T^3
  \ar[ur]_{T \epsilon T}
 } \quad = \quad
\cd{
 T
  \ar[rr]^{\Delta} & &
 T^2
  \ar[rr]^{\id_{T^2}}
  \ar[dr]_{\Delta T} &
 {}
  \dtwocell{d}{\lambda' T} &
 T^2\text. \\ & & &
 T^3
  \ar[ur]_{T \epsilon T}
 }
\]
\end{enumerate}
\end{Defn}

\begin{Defn}
Let $(S, \eta, \mu, \lambda, \rho, \tau)$ be a pseudomonad and $(T, \epsilon,
\Delta, \lambda', \rho', \tau')$ a pseudocomonad on a bicategory $\mathcal B$.
Then a
\defn{pseudo-distributive law} $\delta$ of $T$ over $S$ is given by the following data:
\begin{enumerate}[(PDD1)]
\item A pseudo-natural transformation $\delta \colon TS \Rightarrow ST$;
\item Invertible modifications
\[
\cd[@+1em]{
 T
  \ar[d]_{T\eta}
  \ar[dr]^{\eta T}
  & {} \\
 TS
  \rtwocell[0.35]{ur}{\overline \eta}
  \ar[r]_\delta &
 ST
} \quad \text{and} \quad \cd[@+1em]{
 TS
  \rtwocell[0.35]{dr}{\overline \epsilon}
  \ar[r]^\delta
  \ar[d]_{\epsilon S} &
 ST\text;
  \ar[dl]^{S\epsilon} \\
 S & {}
}\]
\item Invertible modifications
\[
\cd[@+1em]{
 TSS
  \dtwocell{drr}{\overline \mu}
  \ar[r]^{\delta S}
  \ar[d]_{T\mu} &
 STS
  \ar[r]^{S\delta} &
 SST
  \ar[d]^{\mu T} \\
 TS
  \ar[rr]_\delta & &
 ST
} \quad \text{and} \quad \cd[@+1em]{
 TS
  \dtwocell{drr}{\overline \Delta}
  \ar[d]_{\Delta S}
  \ar[rr]^\delta & &
 ST
  \ar[d]^{S \Delta} \\
 TTS
  \ar[r]_{T\delta} &
 TST
  \ar[r]_{\delta T} &
 STT\text,
}
\]
\end{enumerate}
subject to the following axioms
\[\cd{
 TS
  \ar[rr]^{\epsilon S}
  \ar[ddrr]^{\delta} & &
 S \\ &
 {}
  \dtwocell{ur}{\overline \epsilon}
 \\
 T
  \ar[uu]^{T \eta}
  \ar[rr]_{\eta T}
  \dtwocell{ur}{\overline \eta}
 & &
 ST
  \ar[uu]_{S \epsilon}
} = \cd[@-0.6em]{
 TS
  \ar[rr]^{\epsilon S}
  \dtwocell{dr}{\cong} & &
 S \\ &
 \id_{\mathcal B}
  \dtwocell{dr}{\cong}
  \ar[ur]_{\eta} \\
 T
  \ar[uu]^{T \eta}
  \ar[rr]_{\eta T}
  \ar[ur]_{\epsilon}
 & &
 ST
  \ar[uu]_{S \epsilon}
} \eqno{\textrm{(PDA1)}}
\]

\[
\cd[@+1.7em]{
 TSS
  \ar[r]^{\delta S}
  \ar[dr]^{T\mu} &
 STS
  \ar[r]^{S \delta}
  \dtwocell{dr}{\overline \mu} &
 SST
  \ar[d]^{\mu T} \\
 TS
  \dtwocell[0.3]{ur}{T \rho}
  \ar[u]^{T \eta S}
  \ar[r]_{\id_{TS}} &
 TS
  \ar[r]_{\delta} &
 ST
} = \cd[@+1.7em]{
 TSS
  \dtwocell[0.3]{dr}{\overline \eta S}
  \ar[r]^{\delta S} &
 STS
  \ar[r]^{S \delta}
  \dtwocell{d}{\cong} &
 SST
  \ar[d]^{\mu T} \\
 TS
  \ar[ur]_{\eta TS}
  \ar[u]^{T \eta S}
  \ar[r]_{\delta} &
 ST
  \ar[ur]^{\eta ST}
  \ar[r]_{\id_{ST}} &
 ST
  \dtwocell[0.3]{ul}{\rho T}
} \eqno{\textrm{(PDA2)}}
\]

\[
\cd[@+1.7em]{
 TSS
  \ar[r]^{\delta S}
  \ar[dr]^{T\mu} &
 STS
  \ar[r]^{S \delta}
  \dtwocell{dr}{\overline \mu} &
 SST
  \ar[d]^{\mu T} \\
 TS
  \dtwocell[0.3]{ur}{T\lambda}
  \ar[u]^{TS\eta}
  \ar[r]_{\id_{TS}} &
 TS
  \ar[r]_{\delta} &
 ST
} = \cd[@+1.7em]{
 TSS
  \dtwocell{dr}{\cong}
  \ar[r]^{\delta S} &
 STS
  \ar[r]^{S \delta} &
 SST
  \ar[d]^{\mu T} \\
 TS
  \ar[u]^{TS\eta}
  \ar[r]_{\delta} &
 ST
  \ar[u]^{ST\eta}
  \ar[ur]|{S\eta T}
  \ar[r]_{\id_{ST}} &
 ST
  \dtwocell[0.75]{ul}{S \overline \eta}
  \dtwocell[0.3]{ul}{\lambda T}
} \eqno{\textrm{(PDA3)}}
\]

\[
\cd[@+1.7em]{
 TS
  \ar[r]^{\delta}
  \ar[d]_{\Delta S} &
 ST
  \ar[r]^{\id_{ST}}
  \ar[dr]_{S \Delta}
  \dtwocell{dl}{\overline \Delta} &
 ST
  \dtwocell[0.35]{dl}{S \rho'} \\
 TTS
  \ar[r]_{T \delta} &
 TST
  \ar[r]_{\delta T} &
 STT
  \ar[u]_{S \epsilon T}
} = \cd[@+1.7em]{
 TS
  \dtwocell[0.3]{dr}{\rho' S}
  \ar[r]^{\id_{TS}}
  \ar[d]_{\Delta S} &
 TS
  \ar[r]^{\delta}
  \dtwocell{d}{\cong} &
 ST \\
 TTS
  \ar[ur]_{\epsilon TS}
  \ar[r]_{T \delta} &
 TST
  \ar[ur]^{\epsilon ST}
  \ar[r]_{\delta T} &
 STT
  \dtwocell[0.35]{ul}{\overline \epsilon T}
  \ar[u]_{S \epsilon T}
} \eqno{\textrm{(PDA4)}}
\]

\[
\cd[@+1.7em]{
 TS
  \ar[r]^{\delta}
  \ar[d]_{\Delta S} &
 ST
  \ar[r]^{\id_{ST}}
  \ar[dr]_{S \Delta}
  \dtwocell{dl}{\overline\Delta} &
 ST
  \dtwocell[0.35]{dl}{S\lambda'} \\
 TTS
  \ar[r]_{T\delta} &
 TST
  \ar[r]_{\delta T} &
 STT
  \ar[u]_{ST\epsilon}
} = \cd[@+1.7em]{
 TS
  \dtwocell[0.3]{dr}{\lambda' S}
  \ar[r]^{\id_{TS}}
  \ar[d]_{\Delta S} &
 TS
  \dtwocell{dr}{\cong}
  \ar[r]^{\delta} &
 ST \\
 TTS
  \ar[ur]|{T\epsilon S}
  \ar[r]_{T\delta} &
 TST
  \dtwocell[0.35]{ul}{T\overline\epsilon}
  \ar[u]_{TS\epsilon}
  \ar[r]_{\delta T} &
 STT
  \ar[u]_{ST\epsilon}
} \eqno{\textrm{(PDA5)}}
\]

\[
\cd[@R+1em@C+0.3em]{
 & TSSS
  \ar[r]^{\delta SS}
  \ar[dl]_{TS\mu}
  \ar[dr]^{T\mu S} &
 STSS
  \dtwocell{dr}{\overline \mu S}
  \ar[r]^{S\delta S} &
 SSTS
  \ar[r]^{SS\delta}
  \ar[dr]^{\mu TS} &
 SSST
  \dtwocell{d}{\cong}
  \ar[dr]^{\mu ST} \\
 TSS
  \ltwocell{rr}{T \tau}
  \ar[dr]_{T\mu} & &
 TSS
  \ar[rr]_{\delta S}
  \ar[dl]^{T \mu} &
 {}
  \dtwocell{d}{\overline \mu} &
 STS
  \ar[r]_{S \delta} &
 SST
  \ar[dl]^{\mu T} \\ &
 TS
  \ar[rrr]_{\delta} & &
 {}
  \ar@{}[dl]|{\big|\,\big|} &
 ST \\
 & TSSS
  \dtwocell{d}{\cong}
  \ar[r]^{\delta SS}
  \ar[dl]_{TS\mu} &
 STSS
  \dtwocell{dr}{S\overline \mu}
  \ar[dl]_{ST\mu}
  \ar[r]^{S\delta S} &
 SSTS
  \ar[r]^{SS\delta} &
 SSST
  \ar[dl]_{S\mu T}
  \ar[dr]^{\mu ST} \\
 TSS
  \ar[r]_{\delta S}
  \ar[dr]_{T\mu} &
 STS
  \ar[rr]_{S\delta} &
 {}
  \dtwocell{d}{\overline \mu} &
 SST
  \ar[dr]_{\mu T}
  \ltwocell{rr}{\tau T}
  & &
 SST
  \ar[dl]^{\mu T} \\ &
 TS
  \ar[rrr]_{\delta} & & &
 ST
} \eqno{\textrm{(PDA6)}}
\]

\[
\cd[@R+1em@C+0.3em]{
 & TS
  \ar[dl]_{\Delta S}
  \ar[dr]^{\Delta S}
  \ar[rrr]_{\delta} & &
 {}
  \dtwocell{d}{\overline \Delta} &
 ST
  \ar[dr]^{S \Delta} \\
 TTS
  \ltwocell{rr}{\tau' S}
  \ar[dr]_{\Delta TS} & &
 TTS
  \dtwocell{dr}{T\overline \Delta}
  \ar[rr]_{T\delta}
  \ar[dl]^{T\Delta S} & &
 TST
  \dtwocell{d}{\cong}
  \ar[dl]^{TS \Delta}
  \ar[r]_{\delta T} &
 STT
  \ar[dl]^{ST \Delta} \\ &
 TTTS
  \ar[r]_{TT\delta} &
 TTST
  \ar[r]_{T\delta T} &
 TSTT
  \ar@{}[dl]|{\big|\,\big|}
  \ar[r]_{\delta TT} &
 STTT \\
 & TS
  \ar[dl]_{\Delta S}
  \ar[rrr]_{\delta} &
 {}
  \dtwocell{d}{\overline \Delta} & &
 ST
  \ar[dl]_{S \Delta}
  \ar[dr]^{S \Delta} \\
 TTS
  \ar[r]_{T\delta}
  \ar[dr]_{\Delta TS} &
 TST
  \dtwocell{d}{\cong}
  \ar[rr]_{\delta T}
  \ar[dr]_{\Delta ST} &
 {}
  \dtwocell{dr}{\overline \Delta T} &
 STT
  \ltwocell{rr}{S \tau'}
  \ar[dr]_{S \Delta T}
  & &
 STT
  \ar[dl]^{ST \Delta} \\ &
 TTTS
  \ar[r]_{TT\delta} &
 TTST
  \ar[r]_{T\delta T} &
 TSTT
  \ar[r]_{\delta TT} &
 STTT
} \eqno{\textrm{(PDA7)}}
\]

\[\cd[@+0.5em]{
 TSS
  \ar[rrr]^{\epsilon SS}
  \ar[dr]_{\delta S}
  \ar[ddd]_{T\mu} &
 {}
  \dtwocell{d}{\overline \epsilon S} & &
 SS
  \ar[ddd]^{\mu} \\ &
 STS
  \ar[dr]_{S \delta}
  \ar[urr]^{S \epsilon S}
  \dtwocell{rr}{S \overline \epsilon} & & \\ & &
 SST
  \dtwocell{ll}{\overline \mu}
  \ar[d]_{\mu T}
  \ar[uur]_{SS \epsilon}
  \dtwocell{r}{\cong} & \\
 TS
  \ar[rr]_\delta & &
 ST
  \ar[r]_{S \epsilon} &
 S
} \quad = \quad \cd[@+1em]{
 TSS
  \ar[rr]^{\epsilon SS}
  \ar[d]_{T \mu} &
 {}
  \dtwocell[0.4]{dd}{\cong} &
 SS
  \ar[dd]^\mu \\
 TS
  \ar[drr]^{\epsilon S}
  \ar[d]_{\delta} & & \\
 ST
  \ar[rr]_{S \epsilon}
  \dtwocell[0.3]{urr}{\overline \epsilon} & &
 S
} \eqno{\textrm{(PDA8)}}
\]

\[\cd[@+0.5em]{
 T
  \ar[r]^{T \eta}
  \ar[ddd]_{\Delta} &
 TS
  \ar[d]^{\Delta S}
  \ar[rr]^{\delta} & &
 ST
  \ar[ddd]^{S \Delta} \\ &
 TTS
  \dtwocell[0.6]{l}{\cong}
  \dtwocell{rr}{\overline \Delta}
  \ar[dr]^{T \delta} & & \\ & &
 TST
  \ar[dr]^{\delta T}
  \dtwocell{ll}{T \overline \eta}
  \dtwocell{d}{\overline \eta T} \\
 TT
  \ar[uur]^{TT \eta}
  \ar[urr]_{T \eta T}
  \ar[rrr]_{\eta TT} & & &
 STT
} \quad = \quad \cd[@+1em]{
 T
  \ar[rr]^{T \eta}
  \ar[dd]_{\Delta}
  \ar[drr]_{\eta T} & &
 TS
  \ar[d]^{\delta}
  \dtwocell[0.3]{dll}{\overline \eta} \\ & &
 ST
  \ar[d]^{S \Delta} \\
 TT
  \ar[rr]_{\eta TT} &
 {}
  \dtwocell[0.4]{uu}{\cong} &
 STT
} \eqno{\textrm{(PDA9)}}
\]

\[\cd[@!C@!R@-0.4em]{
 & TSS
  \ar[rrr]^{\delta S}
  \ar[ddl]_{T \mu}
  \ltwocell{dddd}{\cong}
  \ar[ddr]^{\Delta SS} &
 {}
  \dtwocell{ddr}{\overline \Delta S} & &
 STS
  \ar[d]^{S \Delta S}
  \ar[rr]^{S \delta} & &
 SST
  \dtwocell{ddl}{S \overline \Delta}
  \ar[ddr]^{SS \Delta} \\ & & & &
 STTS
  \ar[dr]^{ST \delta}
  \dtwocell{dd}{\cong} \\
 \ TS\
  \ar[ddr]_{\Delta S}
  & &
 TTSS
  \ar[r]_{T \delta S}
  \ar[ddl]^{TT\mu} &
 TSTS
  \dtwocell{ddl}{T \overline \mu}
  \ar[ur]^{\delta TS}
  \ar[dr]_{TS \delta} & &
 STST
  \ar[rr]^{S \delta T}
  \dtwocell{ddr}{\overline \mu T} & &
 SSTT
  \ar[ddl]^{\mu TT} \\ & & & &
 TSST
  \ar[ur]_{\delta ST}
  \ar[d]^{T \mu T} \\ &
 TTS
  \ar[rrr]_{T \delta} & & &
 TST
  \ar@{}[dd]|{\big|\,\big|}
  \ar[rr]_{\delta T} & &
 STT \\ \\
 & TSS
  \ar[rrr]^{\delta S}
  \ar[ddl]_{T \mu} & &
 {}
  \dtwocell{dd}{\overline \mu} &
 STS
  \ar[rr]^{S \delta} & &
 SST
  \ar[ddl]_{\mu T}
  \ar[ddr]^{SS \Delta}
  \ltwocell{dddd}{\cong} \\ \\
 TS
  \ar[ddr]_{\Delta S}
  \ar[rrrrr]_\delta & & &
 {}
  \dtwocell{dd}{\overline \Delta} & &
 ST
  \ar[ddr]_{S \Delta} & &
 SSTT
  \ar[ddl]^{\mu TT} \\ \\ &
 TTS
  \ar[rrr]_{T \delta} & & &
 TST
  \ar[rr]_{\delta T} & &
 STT
} \eqno{\textrm{(PDA10)}}\]
\end{Defn}

\bibliography{biblio}

\end{document}